\documentclass[openbib]{svjour3}            
\smartqed  
\usepackage{float}
\usepackage{mathtools}
\usepackage[mathscr]{euscript}
\usepackage{graphicx}
\usepackage{mathptmx}      
\usepackage{chicago-bibstyle}  
\usepackage[misc]{ifsym}
\usepackage{bbding}
\usepackage{appendix}
\usepackage{xcolor}
 \journalname{Computational and Mathematical Biophysics}
\begin{document}

\title{Nonlinear simulation of an elastic tumor-host interface 
}


\author{Min-Jhe Lu$^1$ \and
        Chun Liu$^1$  \and
        Shuwang Li$^1$
}

\authorrunning{Min-Jhe Lu et al.} 

\institute{
\begin{itemize}
    \renewcommand{\labelitemi}{\Letter}
   \item  Shuwang Li \\
            sli@math.iit.edu \\
    \renewcommand{\labelitemi}{$^1$}
    \item Department of Applied Mathematics, 
          Illinois Institute of Technology, Chicago, IL 60616, United States\\
\end{itemize}
}

\date{Received: date / Accepted: date}

\maketitle

\begin{abstract}

We develop a computational method for simulating the nonlinear dynamics of an elastic tumor-host interface. This work is motivated by the recent linear stability analysis of a two-phase tumor model with an elastic membrane interface in 2D \cite{turian2018morphological}. Unlike the classic tumor model with surface tension, the elastic interface condition is numerically challenging due to the 4th order derivative from the Helfrich bending energy.
 Here we are interested in exploring the nonlinear interface dynamics in a sharp interface framework.  We consider a curvature dependent bending rigidity (curvature weakening \cite{he2012modeling}) to investigate  metastasis patterns such as chains or fingers that invade the host environment. We solve the nutrient field and the Stokes flow field using a spectrally accurate boundary integral method, and update the interface using a nonstiff semi-implicit approach. Numerical results suggest curvature weakening promotes the development of branching patterns instead of encapsulated morphologies in a long period of time. For non-weakened bending rigidity, we are able to find self-similar shrinking morphologies based on marginally stable value of the apoptosis rate. 


\end{abstract}

\section{Introduction}\label{intro}
 A malignant tumor usually develops in a sequence of increasingly aggressive stages: carcinogenesis, avascular growth,  angiogenesis and vascular growth \cite{macklin2007nonlinear}.  Avascular growth occurs as tumor cells proliferate and form an \textit{in situ} cancer. Prior to vascularization of the tumor, nutrients (e.g oxygen/glucose) are supplied  through diffusion in the surrounding microenvironment, which limits the size of a tumor spheroid. Morphological instability, however, brings more available nutrients to the tumor by increasing its surface area to volume ratio. In particular, regions where instability first occurs continue to grow at a faster rate than the rest of the tumor tissue. Such diffusional instability is induced by non-uniform cell proliferation and migration according to a heterogeneous distribution of nutrients. The tumor morphology is thus determined by the dominant nutrient levels where proliferation would be favored   \cite{cristini2005morphologic}.

In the past several decades, mathematical models based on fluid mechanics were developed to understand the bio-mechanical properties of the tumor and its metastasis patterns. For example, the original model using Darcy's law (flow through a porous media) \cite{greenspan1976growth,roose2007mathematical} is composed of two parts. One is the concentration of a generic nutrient (e.g. oxygen or glucose)  function $\sigma$ satisfying a reaction diffusion equation $\sigma_t=D \Delta \sigma - \lambda_u \sigma$, where $D$ is the diffusion constant and $\lambda_u$ is the uptake rate; the other part is the internal pressure field $p$ for tumor cell proliferation, which is related to cell velocity $\mathbf{v}$ by the Darcy's law $\mathbf{v}=-\mu \nabla p$, where $\mu$ is the cell mobility. The two is linked by the mass conservation of incompressible tumor cells $\nabla \cdot \mathbf{v} = \lambda_p(\sigma)$, where $\lambda_p(\sigma)=b \sigma-\lambda_A$ is the cell proliferation rate.  $\lambda_M=b \sigma^{\infty}$ and $\lambda_A$ are the rates of mitosis (cell birth) and apoptosis (cell death), respectively. This model is closed by introducing the Laplace-Young condition for internal pressure $(p)_{\partial \Omega}=\gamma \kappa$ (derived from the surface energy of the tumor-host interface), far-field boundary condition for nutrient $(\sigma)_{\partial\Omega}=\sigma^{\infty}$, and the normal velocity of the moving interface $V= -\mu(\nabla p)_{\partial \Omega} \cdot \mathbf{n}$ , where $\gamma$ is the surface tension coefficient and $\kappa$ is the mean curvature. The Darcy flow model, though simple in formulation, captures fundamental features of tumor mechanics and serves as a foundation for developing more sophisticated  models, e.g. bifurcation behavior of the tumor growth by Stokes equation \cite{Hu20071,Hu20072}.

In \cite{pham2018nonlinear}, we introduced a two-phase Stokes model and treated tumor and its host as viscous fluids with different viscosity.  The viscosities reflect the combined properties of cell and extracellular matrix mixtures. The tumor cell population is assumed to be homogeneous and cell proliferation produces a pressure field for tumor growth.   Under the quasi-steady state assumption, this two-phase tumor model consists of a modified Helmholtz equation for nutrient diffusion in tumor and a Stokes equation for tumor dynamics.    In \cite{turian2018morphological}, we extended the two-phase Stokes model by introducing an elastic  tumor-host interface governed by the Helfrich bending energy \cite{Helfrich}.  We derived a modified Laplace-Young condition of the stress jump across the interface for the Stokes equation using an energy variation approach, and performed a linear stability analysis to show how physical parameters such as viscosity, bending rigidity and apoptosis contribute to the morphological instability. Linear results suggest that increased bending rigidity versus mitosis rate contributes to a more stable morphological behavior, and there may exist fingering patterns for increasing tumor viscosity or apoptosis rate.  Comparison with experimental data on glioblastoma spheroids shows good agreement, especially for tumors with high adhesion and low proliferation.

In this paper,  we investigate the nonlinear dynamics of an elastic tumor-host interface. We introduce the bending rigidity coefficient as a function of local mean curvature, referred as curvature weakening model to describe  broken intermolecular bonds and reduced stiffness of the interface  \cite{he2012modeling}. This is also motivated by results from recent studies that changes in stiffness of extracellular matrix may lead to increased mitosis and migration \cite{Mason2012}.   We reformulate the nutrient and Stokes equations as boundary integrals, and develop a sharp interface approach to explore the nonlinear instability of the interface. Note that in the sharp interface framework, we have to compute the 4th order derivative in the interface condition explicitly. After reformulation, the original two-dimension problem is reduced to one-dimensional curve integrals, which can be evaluated using spectrally accurate quadratures. To add more efficiency to the whole algorithm,  a non-stiff interface updating scheme based on the small scale decomposition is implemented  \cite{hou1994removing}.  Our numerical method is spectrally accurate in space and 2nd order accurate in time. Nonlinear simulations show that curvature weakening promotes the development of branching patterns and inhibits encapsulated morphologies. For non-weakened bending rigidity, there exist self-similar shrinking morphologies once the time dependent apoptosis rate (marginally stable value of the apoptosis rate) is applied. Though preliminary, the self-similar idea helps shed light on the strategy for morphological control, as a time dependent apoptosis might be enforced by a well-designed chemo- or radiotherapy.

This paper is organized as follows. In Sect. \ref{sec:2}, we formulate the sharp interface model and non-dimensionalize the resulting PDE systems. In Sect. \ref{sec:3}, we develop BIM formulation and present our numerical method including layer potential evaluations for boundary integrals and small-scale decomposition to remove stiffness. In Sect. \ref{sec:4}, we show numerical results, and then we conclude with Sect. \ref{conclusion}.

\section{Mathematical model}
\label{sec:2}
\begin{figure}[H]
  \centering
  \includegraphics[width=0.5\textwidth]{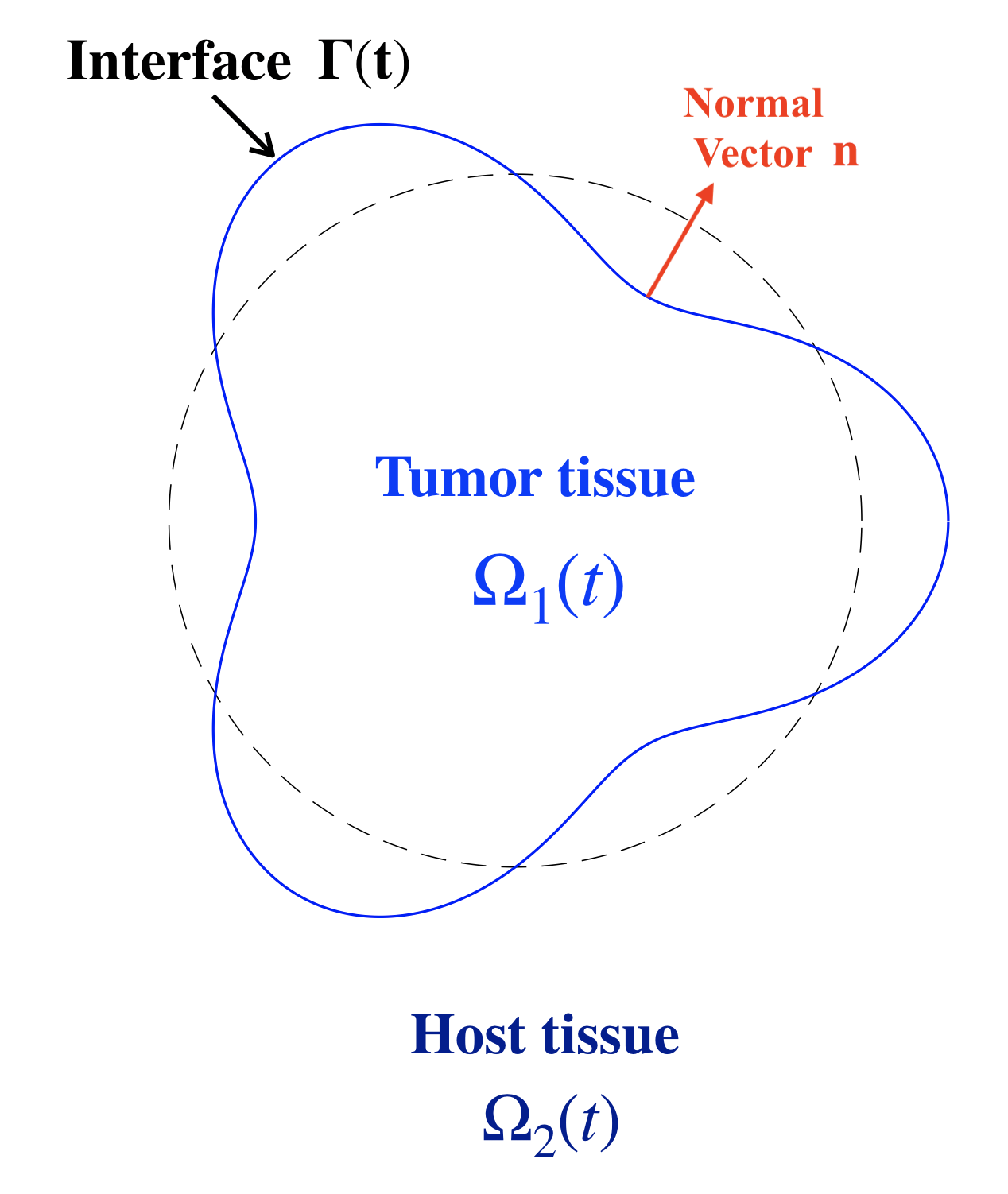}
\caption{Illustration of the computation domain of a three-mode tumor-host interface}
\label{fig:1}     
\end{figure}

We consider an avascular two-dimensional tumor as illustrated in Fig. \ref{fig:1}.  Let $\Omega_1(t)$ be the tumor and $\Omega_2(t)$ be the host tissue. The tumor-host interface $\Gamma(t)$ is considered to be sharp and modeled as an elastic membrane.  At the interface,  a homogeneous elastic bending energy has been widely used to describe the interface dynamics either in a sharp or diffuse interface framework, see e.g.  \cite{Zhou20171,Zhou20172,Promislow2012,DU2004450,Du2005,Wei2010,Dai2013,Sohn2012AxisymmetricMV,Kai2014} among many others---i.e., the elastic bending energy,  $\displaystyle E_{H}=\frac{1}{2} \int_{\Gamma(t)} \nu_{0} \kappa^{2} ds$, where $\nu_0$ is the constant bending rigidity coefficient, $\kappa$ is the local mean curvature, and arc-length $s$ parameterizes the interface. In general,  one may consider a space  dependent energy $\displaystyle E_{W}=\frac{1}{2} \int_{\Gamma(t)} \nu(\kappa) \kappa^{2} ds$, where the bending rigidity $\nu$ is given by a curvature weakening model  \cite{he2012modeling,zhao2016nonlinear},
\begin{equation}
\label{weakening}
\nu(\kappa)=\nu_{0}\left(C e^{-\lambda_{c}^{2} \kappa^{2}}+1-C\right),
\end{equation}
where $0\leq C < 1$ is the rigidity fraction and $\lambda_c$ is the characteristic radius beyond which $\nu(\kappa)$ decays significantly.  When $\kappa$ gets large, parameter $\nu(\kappa)$ approaches to its lower bound $(1-C)\nu_0$. The largest rigidity limit $\nu_0$ can be reached by setting $C = 0$. Performing domain variation in normal direction, i.e. compute variation  \(\frac{\delta E}{\delta \Gamma_{n}} :=\frac{d}{d \epsilon} E(\mathbf{x}+\epsilon \phi \mathbf{n})\), we obtain
\begin{equation}
\label{Evar}
\begin{aligned}
\frac{\delta E_{H}}{\delta \Gamma_{n}}=&-\nu_{0}\left(\frac{1}{2} \kappa^{3}+\kappa_{s s}\right),\\
\frac{\delta E_{W}}{\delta \Gamma_{n}}=&-\left(\frac{1}{2} \nu^{\prime \prime}+2 \nu^{\prime} \kappa+\nu\right) \kappa_{s s} \\ &-\left(\frac{1}{2} \nu^{\prime \prime \prime} \kappa^{2}+3 \nu^{\prime \prime} \kappa+3 \nu^{\prime}\right) \kappa_{s}^{2} -\left(\frac{1}{2} \nu^{\prime} \kappa+\frac{1}{2} \nu\right) \kappa^{3}, \end{aligned}
\end{equation}
where the subscript $s$ denotes a derivative with respect to the arclength parameter $s$ and the prime notation is for a derivative with respect to $\kappa$.
\subsection{The two-phase Stokes tumor model}
\label{sec:2.1}

\paragraph{Nutrient field.}
Similar to the Darcy's model \cite{cristini2003nonlinear}, here the nutrient field in $\Omega_1(t)$ is governed by a reaction diffusion equation:
\begin{equation}
\sigma_t=D \Delta \sigma - \lambda_u \sigma \text{\quad in\ } \Omega_1(t),
\end{equation}
where $D$ and $\lambda_u$ are the diffusion constant and uptake rate, respectively. 
For simplicity, we assume the nutrient concentration \(\sigma\) is constant in \(\Omega_{2}(t)\) and continuous across the interface $\Gamma(t)$:
\begin{equation}
\begin{aligned}
\sigma&=\sigma^{\infty} &\text{in } &\Omega_2(t),\\
\left[\sigma\right]&=0      &\text{on } &\Gamma(t).
\end{aligned}
\end{equation}
\paragraph{Flow field.}
The mass conservation in \(\Omega_{1}(t)\)  and \(\Omega_{2}(t)\) reads:
\begin{equation}\label{massdarcy}
\begin{aligned}
\nabla  \cdot \mathbf{v}
&=\lambda_M\frac{\sigma}{\sigma_\infty}-\lambda_A &\text{in } \Omega_1(t),\\
\nabla  \cdot \mathbf{v}
&=0 &\text{in } \Omega_2(t),
\end{aligned}
\end{equation}
where $\lambda_M$ and $\lambda_A$ are the mitosis and apoptosis rate, respectively. The Stokes equations in both domains are
\begin{equation}
\nabla \cdot \mathbf{T}_{i}=0, \qquad\text{ in } \Omega_i(t), i=1,2,
\end{equation}
where $\displaystyle \mathbf{T}_{i}=\mu_{i}\left(\nabla \mathbf{v}_{i}+\left(\nabla \mathbf{v}_{i}\right)^{T}\right)+\overline{\mu}_{i}\left(\nabla \cdot \mathbf{v}_{i}\right) \mathbf{I}-p_{i} \mathbf{I}$ are stress tensors for the interior tumor \((i=1)\) and exterior host \((i=2)\), \(p_{i}\) are pressures,  and parameters \(\mu_{i}\) and \(\overline{\mu}_{i}\) are  respectively the shear and bulk viscosity coefficients \cite{turian2018morphological,pham2018nonlinear}. The stress jump condition across the interface is given by the bending energy variation in Eq. (\ref{Evar}), 
\begin{equation}
    [\mathbf{T} \mathbf{n}]
    = \frac{\delta E_H}{\delta\Gamma_n}\mathbf{n} \quad \text{or } \quad  \frac{\delta E_W}{\delta \Gamma_n}\mathbf{n}.
\end{equation}
Note that the
stress tensors take into account the rate of strain, dilatation and pressure.
The normal velocity is simply $V=\mathbf{v} \cdot \mathbf{n}$ at the interface. Here we assume the cell velocity $\mathbf{v}$ is continuous across $\Gamma(t)$, i.e. $[\mathbf{v}]=0.$

\subsection{Non-dimensionalization}
Following \cite{turian2018morphological,pham2018nonlinear},  the dimensional variables are scaled by their characteristic values to yield the following non-dimensional parameters:
\begin{equation}
\widetilde{\mathbf{x}}=\frac{\mathbf{x}}{L}, \quad \widetilde{t}=\lambda_{R} t, \quad
\widetilde{\sigma}=\frac{\sigma}{\sigma^{\infty}} , \quad
\widetilde{p}_{i}=\frac{p_i}{\overline{P_1}}, \quad \widetilde{\mathbf{T}_i}=\frac{\mathbf{T}_i}{\overline{T}_{1}}, \quad \widetilde{s}=\frac{s}{L}, \quad \widetilde{\kappa}=L\kappa, \quad i=1,2,
\end{equation}
\noindent 
where $L=\sqrt{\frac{D}{\lambda_u}}$, $\lambda_R ^{-1}=\lambda_M ^{-1}$, and $\sigma^{\infty}$ are the characteristic diffusion length, time, and nutrient concentration scales, respectively. Also,
\(\overline{P}_{1}=\overline{T}_{1}=\mu_{1} \lambda_{M} .\) Since the tumor volume doubling time scale is typically much larger than the diffusion time scale (e.g. days vs. minutes), we assume  $\lambda_{M} \ll \lambda_u$, which leads to a quasi-steady reaction diffusion equation for nutrient field in the tumor tissue.  Dropping all tildes, the nondimensional Stokes-flow system is given by:
\begin{itemize}
\item In the tumor region $\Omega_1(t)$, we have:
\paragraph{Modified Helmholtz equation for nutrient field}
\begin{equation}\label{nutrient}
\Delta {\sigma}= {\sigma}.
\end{equation}
\paragraph{Stokes equation for flow field}
\begin{equation}\label{tumorstokes}
\nabla \cdot  {\mathbf{T}}_1=0,
\end{equation}
where \(\mathbf{T}_{1}=\nabla \mathbf{v}_{1}+\left(\nabla \mathbf{v}_{1}\right)^{T}-\overline{p}_{1} \mathbf{I}\), \(\overline{p}_{1}=p_{1}-\overline{\lambda} \nabla \cdot \mathbf{v}_{1}\) is a modified pressure and $\overline{\lambda}=\frac{\overline{\mu}_{1}}{\mu_{1}}$ is the ratio between two interior viscosities.
\paragraph{Conservation of tumor mass.}
\begin{equation}\label{mass conservation}
 {\nabla} \cdot  {\mathbf{v}}= {\sigma}- \mathcal{A},
\end{equation}
where $\mathcal{A}=\frac{\lambda_A}{\lambda_M}$ represents the relative rate of cell apoptosis to mitosis.

\item In the host tissue region $\Omega_2(t)$, we have a constant nutrient field:
\begin{equation}
     {\sigma}=\sigma|_{\Gamma(t)}=1,
\end{equation}
and divergence free condition for velocity and stress tensor:
\begin{equation}
    \begin{aligned}
    \nabla \cdot \mathbf{v}_2&=0,\\
    \nabla \cdot \mathbf{T}_2&=0,
    \end{aligned}
\end{equation}
where \(\mathbf{T}_{2}=\lambda\left(\nabla \mathbf{v}_{2}+\left(\nabla \mathbf{v}_{2}\right)^{T}\right)-p_{2} \mathbf{I}\) and \(\lambda=\frac{\mu_{2}}{\mu_{1}}\) is the ratio between the exterior host and interior tumor viscosities.

\item On the tumor-host interface $\Gamma(t)$, we have no-jump boundary conditions for nutrient and velocity field:
\begin{equation}\label{nutrientandvelbc}
\begin{aligned}
\left[  {\sigma} \right]&=0,\\
\left[  {\mathbf{v}} \right]&=0,
\end{aligned}
\end{equation}
 and a jump boundary condition for the stress tensor:
\begin{equation}
 {[\mathbf{T n}]}=-\mathscr{S}^{-1} \left(\frac{1}{2}  {\kappa}^{3}+ {\kappa}_{ {s} {s}}\right) \mathbf{n} \text{\quad\quad from \ } E_H,
\end{equation}
\begin{equation}
\begin{aligned}     {[\mathbf{T n}]}=&-\mathscr{S}^{-1}\left(\left(\frac{1}{2} \nu^{\prime \prime}+2 \nu^{\prime}  {\kappa}+\nu\right)  {\kappa}_{ {s}  {s}} +\left(\frac{1}{2} \nu^{\prime \prime \prime}  {\kappa}^{2}+3 \nu^{\prime \prime}  {\kappa}+3 \nu^{\prime}\right)  {\kappa}_{ {s}}^{2}\right. \\
&+ \left.\left(\frac{1}{2} \nu^{\prime}  {\kappa}+\frac{1}{2} \nu\right) {\kappa}^{3} \right)\mathbf{n}
\end{aligned} \text{\quad\quad from \ } E_W,
\end{equation}
where parameter $\displaystyle \mathscr{S}^{-1} 
=\frac{\nu_0}{\mu_1\lambda_M L^3}$ represents the relative strength of bending rigidity.
\end{itemize}

Ideally, one would like to rewrite the Stokes system to the standard one, in which the velocity fields are divergence free in both the tumor and host regions. This is helpful in the design of numerical methods. To do this, we redefine the tumor cell velocity in $\Omega_1(t)$ as
\begin{equation}\label{divfreevel}
    \mathbf{u}_{1}=\mathbf{v}_{1}-\nabla \sigma+\frac{\mathcal{A} \mathbf{x}}{d},
\end{equation}
where $d = 2$ is the spatial dimension. Using the modified Helmholtz equation {\eqref{nutrient} and the identity $\nabla\cdot \mathbf{x}=d$}, equation \eqref{mass conservation} becomes divergence free:
\begin{equation}
\nabla \cdot \mathbf{u}_{1}=0.
\end{equation}
Thus the PDE system in $\Omega_1(t)$ becomes:
\begin{equation}
\label{reformulatedStokes1}
\left\{\begin{array}{l}{\text { Incompressibility } \nabla \cdot \mathbf{u}_{1}=0,} \\ {\text { Stokes equation } \Delta \mathbf{u}_{1}=\nabla \tilde{p}_{1},} \\ {\text { Nutrient equation } \Delta \sigma=\sigma,}\end{array}\right.
\end{equation}
where \(\tilde{p}_{1}\) is the renamed interior pressure \(\tilde{p}_{1}=\overline{p}_{1}-\nabla \cdot \mathbf{v}_{1}- \sigma\). 

{The PDE system in $\Omega_2(t)$ is:
\begin{equation}
\label{reformulatedStokes2}
\left\{\begin{array}{l}{
\text { Incompressibility } \nabla \cdot \mathbf{v}_{2}=0,} \\
{\text { Stokes equation } \lambda\Delta \mathbf{v}_{2}=\nabla p_2,} \\
{\text { Nutrient equation } \sigma=1.}\end{array}\right.
\end{equation}}
Consequently, the boundary conditions can be rewritten as:
\begin{equation}\label{reformulatedbc}
\left\{\begin{array}{l}{\sigma=1} \\
{\mathbf{v}_{2}\left.(\mathbf{x})\right|_{\Gamma{(t)}}-\mathbf{u}_{1}\left.(\mathbf{x})\right|_{\Gamma {(t)}}
=\nabla\left.\sigma\right|_{\Gamma {(t)}}-\frac{\mathcal{A} \left.\mathbf{x}\right|_{\Gamma {(t)}}}{2}} \\
{\mathbf{T}_{2} \mathbf{n}-\mathbf{T}_{1}^{u} \mathbf{n}
=-\mathscr{S}^{-1}f(\kappa) \mathbf{n}+2 \nabla \nabla \sigma \mathbf{n}-2 \sigma \mathbf{n}-\frac{\mathcal{A}}{d}(2-d) \mathbf{n}}
\end{array}\right.
\end{equation}
where \(T_{1}^{u}=\nabla \mathbf{u}_{1}+\left(\nabla \mathbf{u}_{1}\right)^{T}-\tilde{p}_{1}I\), \(f(\kappa)=\frac{1}{2} \kappa^{3}+\kappa_{ss}\) for $E_H$ and \(f(\kappa)=\left(\frac{1}{2} \nu^{\prime \prime}+2 \nu^{\prime} {\kappa}+\nu\right) {\kappa}_{{s} {s}} +\left(\frac{1}{2} \nu^{\prime \prime \prime} {\kappa}^{2}+3 \nu^{\prime \prime} {\kappa}+3 \nu^{\prime}\right) {\kappa}_{{s}}^{2} +\left(\frac{1}{2} \nu^{\prime} {\kappa}+\frac{1}{2} \nu\right) {\kappa}^{3}\) for $E_W$.
The reformulation requires an explicit evaluation of \(\nabla \nabla \sigma \mathbf{n},\) which can be expressed in terms of the normal derivative of \(\sigma\) by Eq. $\eqref{nutrient}$as
\begin{equation}
\begin{aligned} \mathbf{s} \cdot(\nabla \nabla \sigma(s) \mathbf{n}) &=\frac{d}{d s}(\mathbf{n} \cdot \nabla \sigma(s)), \\ \mathbf{n} \cdot(\nabla \nabla \sigma(s) \mathbf{n}) &=1-\kappa \mathbf{n} \cdot \nabla \sigma(s), \end{aligned}
\end{equation}
where $s$ is the arclength representation of the tumor-host interface. Since the exterior velocity field is already divergence-free, it is unnecessary to reformulate the exterior problem. Note that the reformulated velocity field becomes discontinuous across the interface.

 {Notice that although in Eqs. \eqref{reformulatedStokes1}, \eqref{reformulatedStokes2} both the velocity and the nutrient field are governed by time-independent PDEs, the boundary itself and the boundary conditions in Eq. \eqref{reformulatedbc} are time-dependent, which results in a moving boundary problem.}

\subsection{Review of linear stability analysis.}  
\begin{figure}[H]
    \centering
    \includegraphics[width=\textwidth]{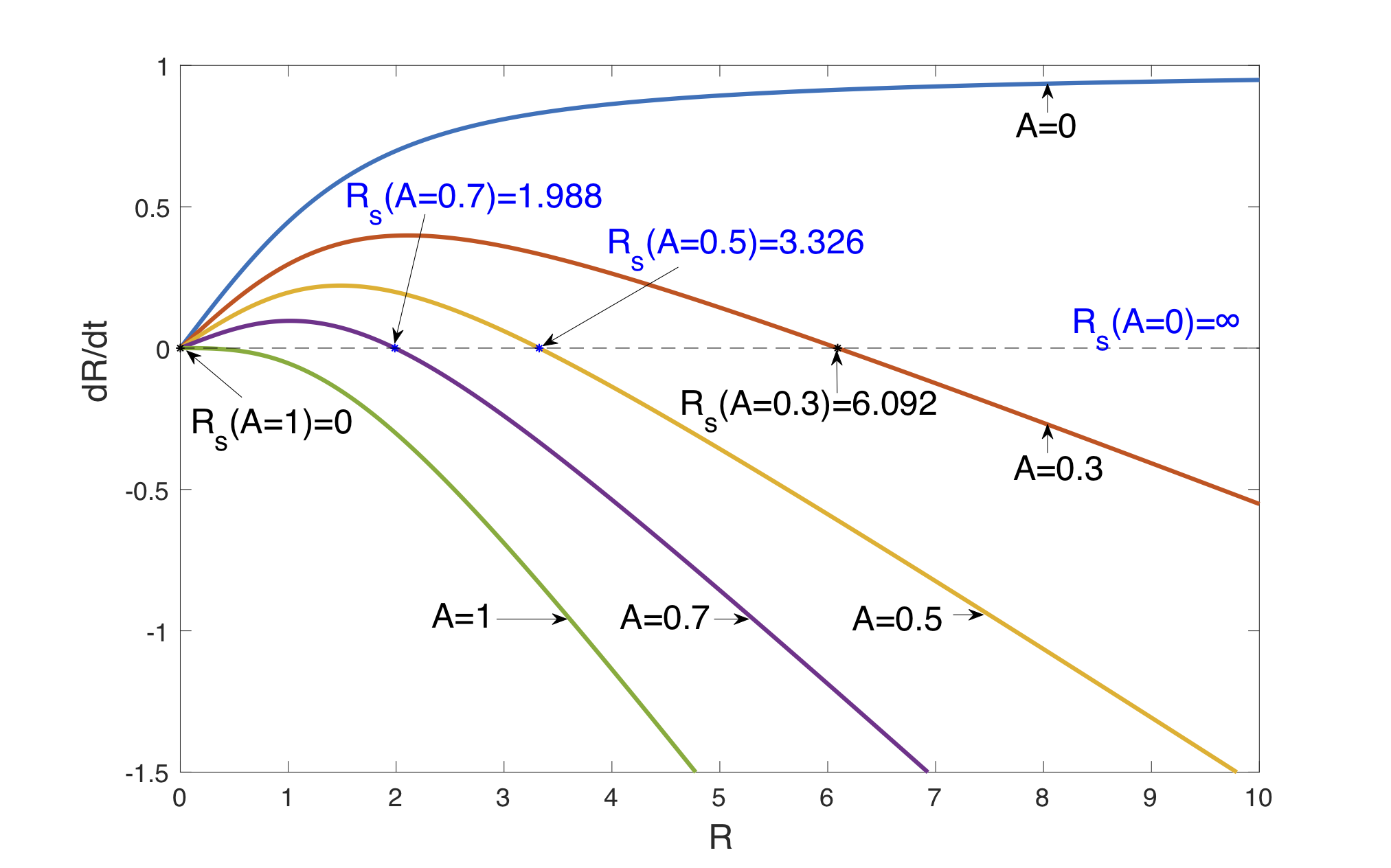}
    \caption{Growth rate from Eq. \eqref{growthrate} for the radially symmetric tumor as a function of $R$.}
    \label{fig:2}
\end{figure}

Though a linear stability analysis could be performed using a weakened bending rigidity, the calculation is very tedious and cumbersome. For brevity, we focus on a constant bending energy case \cite{turian2018morphological} and mainly use our numerical solvers to study the curvature weakening model. For a circular tumor spheroid of radius $R(t)$, the interface evolves as:
\begin{equation}\label{growthrate}
\frac{d R}{d t}=\frac{I_{1}(R)}{I_{0}(R)}-\frac{\mathcal{A} R}{2},
\end{equation}
where $I_0(R)$ and $I_1(R)$ are the modified Bessel functions of the first kind with indices 0 and 1, respectively.
Figure \ref{fig:2} shows the rate of change of a tumor spheroid with respect to the radius $R$ for different $\mathcal{A}$. $R_s(\mathcal{A})$ is the linear steady radius associated with $\mathcal{A}$ satisfying $dR/dt=0$. Large $\mathcal{A}$ causes more cell death and therefore limits the size of the tumor spheroid, whereas $\mathcal{A}=0$ indicates an unbounded growth.
For a slightly perturbed circular interface, $ r(\alpha,t)=R(t)+\epsilon \delta(t) \cos(l \alpha)$, where $\alpha$ is the polar angle, $\delta(t)$ is the time-dependent perturbation, $\epsilon \ll 1$, and
integer $l \ge 2$ is the perturbation mode. The shape perturbation for Stokes-flow with a uniform bending energy evolves as
\begin{equation}\label{shape factor}
\frac{d\left(\frac{\delta}{R}\right)}{d t}=\left(\frac{\delta}{R}\right)\left(\frac{\lambda}{1+\lambda}\mathcal{A}-\frac{l \mathrm{S}^{-1}}{4 R^{3}}\left(l^{2}-\frac{3}{2}\right)+\frac{1}{1+\lambda}\left(1-\frac{I_{1}(R) I_{l+1}(R)}{I_{0}(R) I_{l}(R)}\right)-\frac{2}{R} \frac{I_{1}(R)}{I_{0}(R)}\right),
\end{equation}
where the shape factor $\frac{\delta}{R}$ measures the deviation of the tumor shape from a circle of varying radius, thus describing tumor morphological stability. Note that the normal contribution of the term $\nabla \nabla \sigma$ in the stress jump is focused here. The shape perturbation depends on $\mathcal{A}, \mathscr{S}^{-1}$, $l$ and $\lambda$. Observe that the right hand side of Eq. \eqref{shape factor} increases with increasing $\mathcal{A}$(high cell death) and decreases with increasing $\mathscr{S}^{-1}$ (high membrane rigidity), implying that $\mathcal{A}$ promotes shape instability while $\mathscr{S}^{-1}$ stabilizes it. The morphological stability is mainly determined by the competition between these two important parameters. The parameter $\lambda$ may promote or reduce instability, depending on the values of $l$ and the radius of the tumor $R$. In the result section, we will do a parameter study  on the critical value of the stiffness $\mathscr{S}^{-1}$ and examine the full nonlinear dynamics.

\section{Numerical method}\label{sec:3}

We use the boundary integral methods to solve  (1) the nutrient field in tumor domain; (2) the 2D Stokes equation for the fluid velocity field in both domains. We then update the position of the interface $\Gamma(t)$ by a nonstiff $2{nd}$ order multistep method.  The algorithm presented below is an extension of the approach developed in \cite{pham2018nonlinear} for interfacial flows with surface tension. For completeness, we outline the main ideas here. A rigorous convergence and error analysis of the boundary integral method for the tumor problem can be found in \cite{Wenrui2018}. 
\subsection{Boundary integral formulation}

\paragraph{The nutrient field.}
 Consider a Green's function for the modified Helmholtz equation in $\Omega_1(t)$:
\begin{equation}
\Delta G_{\sigma} - G_{\sigma} = \delta_{\mathbf{x}},
\end{equation}where $G_{\sigma}=G_{\sigma}(\mathbf{x},\mathbf{x}')$, $\mathbf{x}\in \Omega_1(t)$ is the source point, $\mathbf{x}'$ is the field point, and $\delta_{\mathbf{x}}(\mathbf{x},\mathbf{x}')$ is the Dirac delta function. Thus the fundamental solution for Eq. \eqref{nutrient} is the modified Bessel function of the second kind 
\begin{equation}
G_{\sigma}(\mathbf{x},\mathbf{x}')=-\frac{1}{2 \pi} K_0(r),
\end{equation}
where $r=|\mathbf{x}-\mathbf{x}'|$. Using potential theory \cite{kress2013linear}, we define a single-layer potential for the modified Helmholtz equation:
\begin{equation}\label{nutrientsinglelayer}
    (\mathscr{S}_\sigma[\zeta])(\mathbf{x})\coloneqq \int_{\Gamma}G_{\sigma}(\mathbf{x},\mathbf{x}')\zeta(\mathbf{x}')ds',
\end{equation}where $\zeta$ is the layer potential.
The nutrient $\sigma$ can thus be written as a double-layer potential:
\begin{equation}\label{nutrientdoublelayer}
    \sigma(\mathbf{x})=(\mathscr{D}_\sigma[\zeta]) (\mathbf{x})\coloneqq
    \int_{\Gamma} 
    \frac{\partial G_{\sigma}(\mathbf{x},\mathbf{x}')}{\partial \mathbf{n}'} 
    \zeta(\mathbf{x}')ds',
\end{equation}
where $\mathbf{n}'$ is the unit outward normal to $\Gamma(t)$. By the uniform nutrient condition \eqref{reformulatedbc},  we may repose Eq. \eqref{nutrient} as a second-kind Fredholm integral equation with an unknown density $\zeta$ on $\Gamma(t)$:
\begin{equation}
    (-\frac{1}{2}+\mathscr{D}_\sigma)[\zeta]=1.
\end{equation}The term $\frac{\partial \sigma}{\partial \mathbf{n}}$ (normal derivative of double-layer potential) can be computed using
\begin{equation}
    \frac{\partial \sigma}{\partial \mathbf{n}}(s)
    =\frac{d}{d s} \mathscr{S}_\sigma\left[\zeta_{s}\right]-\mathbf{n}(s) \cdot \mathscr{S}_\sigma[\mathbf{n} \zeta].
\end{equation}

\paragraph{The Stokes-flow field.}
Following \cite{pham2018nonlinear}, let \(\mathbf{G}\) be the Stokeslet
and \(\mathbf{T}\) be the tensor stresslet, then the single layer potential at the interface $\Gamma$ is
\begin{equation}
\mathbf{S}[\mathbf{f}]\left(\left.\mathbf{x}\right|_{\Gamma}\right)=\frac{1}{4 \pi} \int_{\Gamma} \mathbf{G}\left(\mathbf{x}^{\prime}-\left.\mathbf{x}\right|_{\Gamma}\right) \mathbf{f}\left(\mathbf{x}^{\prime}\right) d s^{\prime},
\end{equation}
and the double layer potential at $\Gamma$ is
\begin{equation}
{\mathbf{D}}[\mathbf{u}]\left(\left.\mathbf{x}\right|_{\Gamma}\right)=\frac{1}{4 \pi} P . V . \int_{\Gamma} \mathbf{u}\left(\mathbf{x}^{\prime}\right) \mathbf{T}\left(\mathbf{x}^{\prime}-\left.\mathbf{x}\right|_{\Gamma}\right) \mathbf{n}\left(\mathbf{x}^{\prime}\right) d s^{\prime},
\end{equation}
where $P.V.$ indicates the principle value integral.
Assuming that the flow vanishes at the far-field, the boundary integral representation of the 
velocity \(\mathbf{v}_{2}\) approaching the interface from the exterior domain \(\Omega_{2}\) is
\begin{equation}\label{extv}
\mathbf{v}_{2}\left(\left.\mathbf{x}\right|_{\Gamma}\right)=-2 \mathbf{S}\left[\mathbf{f}_{2}\right]\left(\left.\mathbf{x}\right|_{\Gamma}\right)+2 {\mathbf{D}}\left[\mathbf{v}_{2}\right]\left(\left.\mathbf{x}\right|_{\Gamma}\right),
\end{equation}
and the velocity \(\mathbf{u}_{1}\) approaching the interface from the interior domain \(\Omega_{1}\) is
\begin{equation}\label{intv}
\mathbf{u}_{1}\left(\left.\mathbf{x}\right|_{\Gamma}\right)=2 \mathbf{S}\left[\mathbf{f}_{1}^{u}\right]\left(\left.\mathbf{x}\right|_{\Gamma}\right)-2 {\mathbf{D}}\left[\mathbf{u}_{1}\right]\left(\left.\mathbf{x}\right|_{\Gamma}\right),
\end{equation}
where $\mathbf{f}_{i}$ denote the interior and exterior normal stress at the interface \(\left(\lambda^{-1}\mathbf{T}_{2} \mathbf{n}=\mathbf{f}_{2}\right.\)
and \(\mathbf{T}_{1}^{u} \mathbf{n}=\mathbf{f}_{1}^{u} ) .\) Since we do not know \(\mathbf{T}_{i}\) individually, we rewrite Eqs. \eqref{extv} and \eqref{intv} in terms of \(\mathbf{T}_{2}-\mathbf{T}_{1}^{u}\) and make use of Eq. \eqref{reformulatedbc}. To start, we multiply Eq. \eqref{extv} by $\lambda$ and add it to Eq.\eqref{intv} to get
\begin{equation}\label{eq74}
\lambda \mathbf{v}_{2}\left(\left.\mathbf{x}\right|_{\Gamma}\right)+\mathbf{u}_{1}\left(\left.\mathbf{x}\right|_{\Gamma}\right)-2 {\mathbf{D}}\left[\lambda \mathbf{v}_{2}-\mathbf{u}_{1}\right]\left(\left.\mathbf{x}\right|_{\Gamma}\right)=-2 \mathbf{S}\left[ \mathbf{T}_{2}-\mathbf{T}_{1}^{u}\right]\left(\left.\mathbf{x}\right|_{\Gamma}\right),
\end{equation}
where the term \(\mathbf{T}_{2}-\mathbf{T}_{1}^{u}\) is explicitly given in Eq. \eqref{reformulatedbc}. At the interface,  {from Eqs.\eqref{nutrientandvelbc}, \eqref{divfreevel}}, the interior and exterior velocities are related by
\begin{equation}\label{eq75}
\mathbf{v}_{2}\left(\left.\mathbf{x}\right|_{\Gamma}\right)-\mathbf{u}_{1}\left(\left.\mathbf{x}\right|_{\Gamma}\right)=\nabla\left.\sigma\right|_{\Gamma}-\frac{\mathcal{A} \left.\mathbf{x}\right|_{\Gamma}}{2}.
\end{equation}
Putting Eqs. \eqref{eq74} and \eqref{eq75} together, we get
\begin{equation}\label{eq76}
\mathbf{v}_{2}\left(\left.\mathbf{x}\right|_{\Gamma}\right)-2 \frac{\lambda-1}{\lambda+1} {\mathbf{D}}\left[\mathbf{v}_{2}\right]\left(\left.\mathbf{x}\right|_{\Gamma}\right)=\frac{1}{\lambda+1} \mathbf{F},
\end{equation}
where the force term
\begin{equation}
\mathbf{F}=-2 \mathbf{S}\left[\mathbf{T}_{2}-\mathbf{T}_{1}^{u}\right]\left(\left.\mathbf{x}\right|_{\Gamma}\right)+2 {\mathbf{D}}\left[ \left.\nabla \sigma\right|_{\Gamma}-\frac{\mathcal{A} \left.\mathbf{x}\right|_{\Gamma}}{2}\right]\left(\left.\mathbf{x}\right|_{\Gamma}\right)+ \nabla\left.\sigma\right|_{\Gamma}-\frac{\mathcal{A} \left.\mathbf{x}\right|_{\Gamma}}{2}.
\end{equation}
Note that in 2D, \(\mathbf{v}_{2}=\left(v_{1}, v_{2}\right), \mathbf{n}=\left(n_{1}, n_{2}\right),\) and \(\mathbf{F}=\left(F_{1}, F_{2}\right) .\) Using the formulas of
the double and single layer potentials \cite{pozrikidis1992boundary}, equation \eqref{eq76} can be explicitly rewritten as
\begin{equation}\label{stokesbim}
v_{j}\left(\left.\mathbf{x}\right|_{\Gamma}\right)-2 \frac{\lambda-1}{\lambda+1} \frac{1}{4 \pi} \int_{\Gamma} v_{i}\left(\mathbf{x}^{\prime}\right) T_{i j k}\left(\mathbf{x}^{\prime}, \left.\mathbf{x}\right|_{\Gamma}\right) n_{k}\left(\mathbf{x}^{\prime}\right) d s^{\prime}=\frac{1}{\lambda+1} F_{j}\left(\
\left.\mathbf{x}\right|_{\Gamma}\right), \quad j=1,2
\end{equation}
where
\begin{equation}
\begin{aligned} F_{j}\left(\left.\mathbf{x}\right|_{\Gamma}\right)=&-2 \frac{1}{4 \pi} \int_{\Gamma} f_{i}\left(\mathbf{x}^{\prime}\right) G_{i j}\left(\mathbf{x}^{\prime}, \left.\mathbf{x}\right|_{\Gamma}\right) d s^{\prime} \\
&+2 \frac{1}{4 \pi} \int_{\Gamma} h_{i}\left(\mathbf{x}^{\prime}\right) T_{i j k}\left(\mathbf{x}^{\prime}, \left.\mathbf{x}\right|_{\Gamma}\right) n_{k}\left(\mathbf{x}^{\prime}\right) d s^{\prime}+h_{j}\left(\left.\mathbf{x}\right|_{\Gamma}\right), \end{aligned}
\end{equation}
 \(\mathbf{h}\left(\left.\mathbf{x}\right|_{\Gamma}\right)=\nabla \sigma\left(\left.\mathbf{x}\right|_{\Gamma}\right)-\frac{\mathcal{A} \left.\mathbf{x}\right|_{\Gamma}}{2}, G_{i j}=\sum_{i}^{d}\left(-\delta_{i j} \ln r+\frac{\hat{\mathbf{x}}_{i} \hat{\mathbf{x}}_{j}}{r^{2}}\right)\) and \(T_{i j k}=\)
\(\sum_{i, k}^{d}\left(-4 \frac{\hat{\mathbf{x}}_{i} \hat{\mathbf{x}}_{j}\hat{\mathbf{x}}_{k}}{r^{2} }\right) \) with \(r=|\hat{\mathbf{x}}|\) and \(\hat{\mathbf{x}}=\mathbf{x}^{\prime}(s)-\left.\mathbf{x}\right|_{\Gamma}(s) .\) Hence, the explicit forms of the single layer and double layer potentials are
\begin{equation}\label{eq80}
\int_{\Gamma} f_{i}\left(\mathbf{x}^{\prime}\right) G_{i j}\left(\mathbf{x}^{\prime}, \left.\mathbf{x}\right|_{\Gamma}\right) d s^{\prime}=
\left\{\begin{array}{l} \int_{\Gamma}\left(-f_{1}^{\prime} \log r+f_{1}^{\prime} \frac{\hat{x}_{1}^{2}}{r^{2}}+f_{2}^{\prime} \frac{\hat{x}_{1} \hat{x}_{2}}{r^{2}}\right) d s^{\prime} \quad j=1, \\ \int_{\Gamma}\left(-f_{2}^{\prime} \log r+f_{2}'\frac{\hat{x}_{2}^{2}}{r^2}+f_{1}^{\prime} \frac{\hat{x}_{1} \hat{x}_{2}}{r^{2}}\right) d s^{\prime} \quad j=2, \end{array}\right.
\end{equation}
and
\begin{equation}
\int_{\Gamma} v_{i}\left(\mathbf{x}^{\prime}\right) T_{i j k}\left(\mathbf{x}^{\prime}, \left.\mathbf{x}\right|_{\Gamma}\right) n_{k}\left(\left.\mathbf{x}\right|_{\Gamma}\right) d s^{\prime}
=
\left\{\begin{array}{l} 
\int_{\Gamma} \frac{-4}{r^{4}}\left(v_{1}^{\prime} \hat{x}_{1}^{3} n_{1}+v_{1}^{\prime} \hat{x}_{1}^{2} \hat{x}_{2} n_{2}+v_{2}^{\prime} \hat{x}_{1}^{2} \hat{x}_{2} n_{1}+v_{2}^{\prime} \hat{x}_{1} \hat{x}_{2}^{2} n_{2}\right) d s^{\prime} \quad j=1, \\
\int_{\Gamma} \frac{-4}{r^{4}}\left(v_{1}^{\prime} \hat{x}_{1}^{2} \hat{x}_{2} n_{1}+v_{1}^{\prime} \hat{x}_{1} \hat{x}_{2}^{2} n_{2}+v_{2}^{\prime} \hat{x}_{1} \hat{x}_{2}^{2} n_{1}+v_{2}^{\prime} \hat{x}_{2}^{3} n_{2}\right) d s^{\prime} \quad j=2,
\end{array}\right.
\end{equation}
where \(\hat{x}_{1}=x(s(\alpha))-x\left(s\left(\alpha^{\prime}\right)\right), \hat{x}_{2}=y(s(\alpha))-y\left(s\left(\alpha^{\prime}\right)\right), v_{i}^{\prime}=v_{i}\left(\mathbf{x}\left(s\left(\alpha^{\prime}\right)\right)\right)\) and
\(n_{i}=n_{i}(\mathbf{x}(s(\alpha)))\).
In Eq. \eqref{eq80}, the only singularity in the integrand comes from the logarithmic
kernel. This can be analyzed in the following subsection.

\subsection{The evaluation of the boundary integrals \cite{li2011boundary}}
With the integral formulation above, we assume the interface  $\Gamma$ is analytic and given by $\big\{\mathbf{x}(\alpha,t)=(x(\alpha,t),y(\alpha,t): 0\leq \alpha \leq 2 \pi \big\}$, where $\mathbf{x}$ is $2 \pi$-periodic in the parametrization $\alpha$. The unit tangent and normal (outward) vectors can be calculated as $\mathbf{s}=(x_\alpha,y_\alpha)/s_\alpha$, $\mathbf{n}=(y_\alpha,-x_\alpha)/s_\alpha$, where the local variation of the arclength $s_\alpha=\sqrt{x_\alpha^2+y_\alpha^2}$. Subscripts refer to partial differentiation.
We track the interface $\Gamma$ by introducing $N$ marker points to discretize the planar curves, parametrized by $\alpha_j=jh$, $h=\frac{2\pi}{N}$, where $N$ is a power of $2$.
\paragraph{Computation of the single-layer potential.\\}
In Eqs. \eqref{nutrientsinglelayer} and \eqref{eq80}, the single-layer potential type integrals contain the Green functions with the logarithmic singularity at $r=0$. They can be rewritten as the following form under the parametrization $\alpha$:
\begin{equation}\label{singlelayeralpha}
    \int_{\Gamma} 
    \Phi(\alpha,\alpha')\phi(\alpha')s_{\alpha}(\alpha')d\alpha',
\end{equation}where $\Phi$ are the Green functions $G$ or $G_\sigma$, $\phi$ is the layer density $\eta$ or $\zeta$. We  decompose the Green functions as:
\begin{equation}\label{decomposegreenp}
    G(\alpha,\alpha')=
    -\frac{1}{2 \pi} \ln r
    =-\frac{1}{2\pi}\left(
    \ln{2 \left| \sin{\frac{\alpha-\alpha'}{2}}\right|}
    +\left[\ln r-\ln{2 \left| \sin{\frac{\alpha-\alpha'}{2}}\right|}\right]
    \right),
\end{equation}

\begin{equation}\label{decomposegreensigma}
    G_\sigma(\alpha,\alpha')=
    -\frac{1}{2 \pi} K_{0}(r)=
    -\frac{1}{2 \pi}
    \left(I_{0}(r)\ln{2 \left| \sin{\frac{\alpha-\alpha'}{2}}\right|}
    +\left[K_{0}(r)
    -I_{0}(r)\ln{2 \left| \sin{\frac{\alpha-\alpha'}{2}}\right|}
    \right]\right),
\end{equation}where $I_{0}$ is a modified Bessel function of the first kind, $r=|\mathbf{x}(\alpha)-\mathbf{x}'(\alpha ')|$. The square bracket on the right-hand side of Eqs.\eqref{decomposegreenp} and \eqref{decomposegreensigma} has removable singularity at $\alpha=\alpha'$, since $r=
s_{\alpha}\left|\alpha-\alpha'\right|
\sqrt{1+\mathscr{O}(\alpha-\alpha')}
=s_{\alpha}\left|\alpha-\alpha'\right|
(1+\mathscr{O}(\alpha-\alpha'))$ for $\alpha \approx \alpha'$, where $\mathscr{O(\alpha-\alpha')}$ denotes a smooth function that vanishes as $\alpha\rightarrow\alpha'$, and since $K_0$ has the expansion $$K_{0}(z)=-\left(\log \frac{z}{2}+C\right) I_{0}(z)+\Sigma_{n=1}^{\infty} \frac{\psi(n)}{(n !)^{2}}\left(\frac{z}{2}\right)^{2 n}.$$ As a result, for an analytic and $2\pi$-periodic function $f(\alpha,\alpha')$, a standard trapezoidal rule or alternating point rule can be implemented to evaluate the integral 
\begin{equation}
    \int_{0}^{2\pi}
f(\alpha,\alpha')
\ln{\frac{r}{2 \left| \sin{\frac{\alpha-\alpha'}{2}}\right|}}
d\alpha'
\end{equation} and achieve spectral accuracy. The first term on the right-hand side of Eqs.\eqref{decomposegreenp} and \eqref{decomposegreensigma} is still singular and can be evaluated using the following spectrally accurate quadrature \cite{kress1995numerical}:
\begin{equation}\label{logsingularity}
    \int_{0}^{2\pi}f(\alpha_i,\alpha')
\ln{2 \left| \sin{\frac{\alpha_i-\alpha'}{2}}\right|}
d\alpha'\approx
\Sigma_{j=0}^{2m-1}q_{\left|j-i\right|}f(\alpha_i,\alpha_j),
\end{equation}where $m=\frac{N}{2}$, $\alpha_i=\frac{\pi i}{m}$ for $i=0,1,...,2m-1$, and weight coefficients
\begin{equation}\label{Kressquadrature}
    q_j=-\frac{\pi}{m}\Sigma_{k=1}^{m-1}\frac{1}{k}\cos{\frac{kj\pi}{m}}-\frac{(-1)^j\pi}{2m^2} , \text{for } j=0,1,...,2m-1.
\end{equation}The derivative $\frac{d}{d\alpha}$ in Eq. \eqref{singlelayeralpha} is approximated using Fast-Fourier-Transform spectral derivatives thus maintaining spectral accuracy.
\paragraph{Computation of the double-layer potential.\\}
In Eq. \eqref{nutrientdoublelayer}, the double-layer potential type integrals contain the Green functions with logarithmic singularity at $r=0$. It can be rewritten as the following form under the parametrization $\alpha$:
\begin{equation}
    \int_{\Gamma} 
    \frac{\partial \Phi(\alpha,\alpha')}{\partial \mathbf{n}(\alpha')}\phi(\alpha')s_{\alpha}(\alpha')d\alpha',
\end{equation}where $\Phi$ stands for the Green function $G_\sigma$ and $\phi$ is the layer density $\zeta$.

Since $\frac{\partial G_\sigma}{\partial \mathbf{n}}$ has logarithmic singularity, we decompose it as below:
\begin{equation}
    \frac{\partial G_{\sigma}(\alpha,\alpha')}{\partial \mathbf{n}(\alpha')}=
    -h(\alpha,\alpha') K_{1}(\mu_{i} r)=
    g_1(\alpha,\alpha')\ln{2 \left| \sin{\frac{\alpha-\alpha'}{2}}\right|}
    +g_2(\alpha,\alpha'),
\end{equation}where $g_1(\alpha,\alpha')$ and $g_2(\alpha,\alpha')$ are analytic and $2\pi$-periodic functions with
\begin{equation}
    g_1(\alpha,\alpha')=-h(\alpha,\alpha')I_1(\mu_i r),
\end{equation}
\begin{equation}\label{g2}
    g_2(\alpha,\alpha')=-h(\alpha,\alpha')
    \left[K_{1}(\mu_{i}r)
    -I_{1}(\mu_{i}r)\ln{2 \left| \sin{\frac{\alpha-\alpha'}{2}}\right|}
    \right].
\end{equation}We have used the fact $\frac{d}{dr}K_0(r)=-K_1(r)$. Since $K_1$ has the expansion $$K_{1}(z)=\frac{1}{z}+\left(\log \frac{z}{2}+C\right) I_{1}(z)-\frac{1}{2} \sum_{n=0}^{\infty} \frac{\psi(n+1)+\psi(n)}{n !(n+1) !}\left(\frac{z}{2}\right)^{2 n+1},$$ the square bracket on the right-hand side of Eq. $\eqref{g2}$ also has a removable singularity at $\alpha=\alpha'$, thus the integral involving $g_2(\alpha,\alpha')$ can be evaluated by a standard trapezoidal rule or alternating point rule. Note that 
\begin{equation}
    g_2(\alpha,\alpha)=-\frac{h(\alpha,\alpha)}{ r}=-\frac{1}{4\pi}
    \frac{x_\alpha y_{\alpha\alpha}-x_{\alpha\alpha}y_\alpha}{x_\alpha^2+y_\alpha^2}.
\end{equation}
The first term on the right-hand side of Eq. $\eqref{g2}$ is still singular and can be evaluated through the quadrature given in Eqs. \eqref{logsingularity} and \eqref{Kressquadrature}.\\
To summarize, using $\text{Nystr\"om}$ discretization with the Kress quadrature rule \cite{hao2014high} described above, we discretize the boundary integral equations for the nutrient and Stokes fields into two dense linear systems with unknowns as the layer densities $\zeta$ and $\eta$ on $\Gamma(t)$, which can be solved using an iterative solver, e.g., GMRES \cite{saad1986gmres}.
\subsection{The evolution of the interface}
As indicated in \cite{hou1994removing}, the curvature driven motion introduces high-order derivatives, both non-local and non-linear, into the dynamics through the Laplace-Young condition at the interface. Explicit time integration methods suffer from severe numerical stability constraints and implicit methods are difficult to apply since the stiffness enters non-linearly. Hou et al. resolves these difficulties by adopting the $\theta-L$ formulation and the small-scale decomposition (SSD) which we will follow in this paper.
\paragraph{$\theta-L$ formulation.\\}
This description makes the application of an implicit method straightforward and may circumvent the problem of point clustering. Consider a point $\mathbf{x}(\alpha,t)=(x(\alpha,t),y(\alpha,t))\in \Gamma(t)$. Denote the normal and tangent velocity by $V(\alpha,t)={\bf u}\cdot \mathbf{n}$ and $T(\alpha,t)={\bf u}\cdot \mathbf{s}$ respectively, where ${\bf u}=\mathbf{x}_t=V {\mathbf{n}}+T {\mathbf{s}}$ describes the motion of $\Gamma(t)$. The tangent angle that the planar curve $\Gamma(t)$ forms with the horizontal $x$-axis, called $\theta$, satisfies $\theta=\tan^{-1}{\frac{y_\alpha}{x_\alpha}}$. The unit tangent and normal vectors become ${\mathbf{s}}=(\cos{\theta}, \sin{\theta})$ and ${\mathbf{n}}=(\sin{\theta},-\cos{\theta})$. The length of one period of the curve is $L(t)=\int_0^{2\pi}s_\alpha d\alpha$. Differentiating these two equations of $\theta,s_\alpha$ in time, we obtain the following evolution equations:
\begin{equation}\label{evolvetheta}
    \theta_t=\kappa T - V_s=\frac{1}{s_\alpha}(\theta_\alpha T- V_\alpha),
\end{equation}
\begin{equation}\label{evolvearclength}
    s_{\alpha t}=(T_s+\kappa V)s_\alpha=T_\alpha+\theta_\alpha V,
\end{equation}
 {where the curvature is evaluated by $\kappa=\theta_s=\frac{\theta_\alpha}{s_\alpha}.$ }
Instead of using the $(x,y)$ coordinates, we are able to repose the equation of motion in terms of  dynamical variables $(L,\theta)$. 
 To gain more efficiency and accuracy, one may choose a tangent velocity $T$ (independent of the morphology of the interface) such that the marker points are equally spaced in arclength to prevent point clustering:
\begin{equation}
    T(\alpha,t)= \frac{\alpha}{2\pi}\int_0^{2\pi}\theta_{\alpha'}V' d\alpha'-\int_0^\alpha \theta_{\alpha'}V' d\alpha'.
\end{equation}
It follows that $s_\alpha$ is independent of $\alpha$ and thus is everywhere equal to its mean:
\begin{equation}
    s_\alpha=\frac{1}{2\pi}\int_0^{2\pi}s_\alpha(\alpha,t)d\alpha=\frac{L(t)}{2\pi}.
\end{equation}
The procedure for obtaining the initial equal arclength parametrization is presented in Appendix B of \cite{baker1990connection}. The
idea is to solve the nonlinear equation 
\begin{equation}
    \int_0^{\alpha_j} s_{\beta}d\beta=\frac{j}{N}L
\end{equation}for $\alpha_j$ using Newton's method and evaluate the equal arclength marker points $\mathbf{x}(\alpha_j)$ by interpolation in Fourier space. 
We may recover the interface by simply integrating
\begin{equation}\label{xalpha}
    \mathbf{x}_\alpha=\mathbf{x}_s s_\alpha=\frac{L(t)}{2\pi}(\cos{\theta(\alpha,t)},\sin{\theta(\alpha,t)}).
\end{equation}
\paragraph{Small scale decomposition (SSD).\\}
The idea of small scale decomposition (SSD) is to extract the dominant part of the equations at small spatial scales \cite{hou1994removing}. To remove the stiffness, we use SSD in our problem and develop an explicit, non-stiff time integration algorithm.  Through the analysis of the single-layer and double-layer terms, the only singularity in the integrands comes from the logarithmic kernel. Following \cite{hou1994removing} and noticing the curvature terms in the  {stress-jump} condition  {in Eq. \eqref{reformulatedbc} and Eq. \eqref{stokesbim}}, one can show that at small spatial scales  {\cite{Sohn2012AxisymmetricMV}},
\begin{equation}\label{decomposetheta}
\begin{aligned}
    V(\alpha,t) 
    &\sim \frac{1}{s_\alpha^2} \mathcal{H}[\theta_{\alpha\alpha}], 
\end{aligned}
\end{equation}
where $\mathcal{H}(\xi)=\frac{1}{2\pi}\int_0^{2\pi}\xi'\cot{\frac{\alpha-\alpha'}{2}}d\alpha'$ is the Hilbert transform for a $2\pi$-periodic function $\xi$.\\ We rewrite Eq. \eqref{evolvetheta},
\begin{equation}
\label{evolvethetaN}
\begin{aligned}
    \theta_t&=\frac{1}{s_\alpha^3} \mathcal{H}[\theta_{\alpha\alpha\alpha}]+N(\alpha,t),
\end{aligned}
\end{equation}
where the Hilbert transform term is the dominating high-order term at small spatial scales, and $ \displaystyle N= (\kappa T-V_s)-\frac{1}{s_\alpha^3} \mathcal{H}[\theta_{\alpha\alpha\alpha}]$ contains all other lower-order terms in the equation of motion. This splitting reveals that an explicit time-stepping method has the high-order constraint$\displaystyle \left ( \frac{h}{s_\alpha} \right)^{3},$ where $\Delta t$ and $h$ are the time-step and spatial grid size, respectively. The efficiency has been demonstrated numerically in the seminal work \cite{hou1994removing} and later in \cite{ShuwangJCP,Meng2017} for a Hele-Shaw problem.  For the tumor growth problem, the semi-implicit time-stepping scheme (see Eq. \eqref{decomposetheta}) requires $\Delta t = O(h)$ instead of explicit schemes which would require $\Delta t = O(h^3)$.  In section \ref{sec:4}, we show numerical examples using $N=2048$. In the simulation, we could use $\Delta t$ as $\Delta t=1.0 \times 10^{-2}$ for stability instead of $\Delta t<10^{-6}$ for an explicit scheme.
\subsection{Semi-implicit time-stepping scheme}
Taking the Fourier transform of Eq. \eqref{evolvethetaN}, we get 
\begin{equation}\label{fouriertheta}
{\hat \theta}_t=-\left(\frac{|k|}{s_\alpha}\right)^3 {\hat \theta}(k,t) +{\hat N}(k,t).
\end{equation}
In Fourier space, we solve Eq. \eqref{fouriertheta} using a second order accurate linear propagator method in the Adams-Bashforth form and then apply the inverse Fourier transform to recover $\theta$.  Specifically,  we discretize Eq. \eqref{fouriertheta} as
\begin{equation}\label{discevovetheta}
\quad\quad\;{\hat \theta}^{n+1}(k)=e_k(t_n,t_{n+1}){\hat \theta}^{n}(k)+\frac{\Delta t}{2}(3e_k(t_n,t_{n+1}){\hat N}^{n}(k)-e_k(t_{n-1},t_{n+1}){\hat N}^{n-1}(k),
\end{equation}
where  the superscript $n$ denotes the numerical solutions at $t=t_n$ and the integrating factor 
\begin{equation}
e_k(t_1,t_2)=\exp\left (-{|k|^3}\int_{t_1}^{t_2}\frac{dt}{s_\alpha^3(t)}\right ).
\end{equation}
Note that by setting the integrating factors in Eq. \eqref{discevovetheta} to $1$, we recover the classical Adams-Bashforth explicit time-stepping method.
The integrating factors in Eq. \eqref{discevovetheta} can be evaluated simply using trapezoidal rule,
\begin{eqnarray}\label{eq73}
\int_{t_n}^{t_{n+1}}\frac{dt}{s_\alpha^3(t)} &\approx& \frac{\Delta t}{2} \left (\frac{1}{(s_\alpha^n)^3}+\frac{1}{(s_\alpha^{n+1})^3} \right ) \nonumber, \\ 
\int_{t_{n-1}}^{t_{n+1}}\frac{dt}{s_\alpha^3(t)} &\approx& {\Delta t} \left (\frac{1}{2(s_\alpha^{n-1})^3}+\frac{1}{(s_\alpha^{n})^3}+\frac{1}{2(s_\alpha^{n+1})^3} \right ).
\end{eqnarray}
To compute the arclength $s_\alpha$, equation \eqref{evolvearclength} is discretized using the explicit 2nd-order Adams-Bashforth method,
\begin{equation}\label{AB2}
s_\alpha^{n+1}=s_\alpha^n+\frac{\Delta t}{2}(3M^n-M^{n-1}),
\end{equation}
where $M$ is calculated using $M=\frac{1}{2\pi}\int_0^{2\pi}V(\alpha,t)\theta_\alpha d\alpha$.

 Note that the second order linear propagator and Adams-Bashforth methods are multi-step methods and require two previous time steps. The first time step is realized using an explicit Euler method for $s_\alpha^1$ and a first order linear propagator of a similar form for $\hat{\theta}^1$.
 
 To reconstruct the tumor-host interface $(x(\alpha,t_{n+1}),y(\alpha,t_{n+1}))$ from the updated $\theta^{n+1}(\alpha)$ and $s_\alpha^{n+1}$, we first update a reference point $(x(0,t_{n+1}),y(0,t_{n+1}))$ using a second-order explicit Adams-Bashforth method to discretize the equation of motion $\mathbf{x}_t=V\hat{\mathbf{n}}$ (with the tangential part dropped since it does not change the morphology)
 \begin{equation}
     (x(0,t_{n+1}),y(0,t_{n+1}))=(x(0,t_{n}),y(0,t_{n}))+\frac{\Delta t}{2} \left( 3V(0,t_n)\hat{\mathbf{n}}(0,t_{n})-V(0,t_{n-1})\hat{\mathbf{n}}(0,t_{n-1}) \right).
 \end{equation}
 Once we update the reference point, we obtain the configuration of the interface from $\theta^{n+1}(\alpha)$ and $s_\alpha^{n+1}$ by integrating Eq. \eqref{xalpha} following \cite{hou1994removing}
 \begin{eqnarray}
     x(\alpha,t_{n+1})&=&x(0,t_{n+1})+s_\alpha^{n+1}\left( \int_0^{\alpha} \cos(\theta^{n+1}(\alpha'))d\alpha'
     -\frac{\alpha}{2\pi}\int_0^{2\pi}\cos(\theta^{n+1}(\alpha'))d\alpha'\right),\nonumber\\
     y(\alpha,t_{n+1})&=&y(0,t_{n+1})+s_\alpha^{n+1}\left( \int_0^{\alpha} \sin(\theta^{n+1}(\alpha'))d\alpha'
     -\frac{\alpha}{2\pi}\int_0^{2\pi}\sin(\theta^{n+1}(\alpha'))d\alpha'\right),
 \end{eqnarray}where the indefinite integration is performed using the discrete Fourier transform.
 
 We use a 25th order Fourier filter to damp the highest nonphysical mode and suppress the  aliasing error \cite{hou1994removing}. We also use Krasny filtering to prevent the accumulation of round-off errors during the computation \cite{krasny1986study}.

\section{Results}
\label{sec:4}
We have computed a number of different cases which illustrate and expand upon the linear stability analysis in \cite{turian2018morphological}. First, we examine stability of a perturbed interface through the full nonlinear simulations. We then take into account the curvature weakening of the bending rigidity and examine its effects on the pattern formation. Finally, we examine the existence of possible self-similar solution.

The correctness of implementations of boundary integral methods for both Stokes flow and the modified Helmholtz equation was checked in a number of ways. This includes growing/shrinking of a circular interface to the steady radius predicted by Eq. \eqref{growthrate}. Agreement of the nonlinear evolution of perturbations was also verified with linear solution for a slightly perturbed  circular interface in Eq. \eqref{shape factor}. This is assessed by comparing the corresponding linear and nonlinear shape factors. The linear shape factor is calculated by solving Eqs. \eqref{growthrate} and  \eqref{shape factor}, and the nonlinear shape factor is calculated numerically using
\begin{equation}
    \frac{\delta}{R}=\max_{j}(|\frac{\mathbf{x}_j}{R}|^2-1)^{1/2}, j=1,...,N,
\end{equation}
where $\mathbf{x}_j$ denotes the discrete points that describe the tumor/host interface and $R$ denotes the effective radius of the tumor, which is the radius of a circle with the same area as the tumor.  

\subsection{Growth or shrinkage through marginally stable curve $S_M^{-1}$ ($\mathcal{A}=0.5$)}
A marginally stable (or critical) value of the rigidity parameter $S^{-1}_M(l, \mathcal{A}, R, \lambda)$ is obtained by setting the time derivative of $\frac{\delta}{R}$ in Eq. \eqref{shape factor} to zero and thus separates stable $(\mathscr{S}^{-1}>S^{-1}_M)$ from unstable regime $(\mathscr{S}^{-1}<S^{-1}_M)$. Recall that $\mathscr{S}^{-1}$ is proportional to membrane rigidity. In Fig. \ref{fig:marginalSinv} [a], we illustrate this behavior by plotting $\mathscr{S}^{-1}$ with $\mathcal{A}=0.5$, mode $l=3$ for various viscosity ratios $\lambda=0.5,1.5,2.5$ against $R$. When $\mathcal{A}=0.5$, the steady radius is $R_s\approx3.326$. We consider the dynamics of tumors that may grow or shrink depending upon their initial radius. We take two membrane rigidity parameters $\mathscr{S}^{-1} = 0.001$ and $\mathscr{S}^{-1} = 2$, two initial tumor radii $R(0) = 1.988$ and $R(0) = 4.5$, and vary the viscosity ratio $\lambda$. 

Specifically, we consider in Fig. \ref{fig:marginalSinv}  [a] evolution from the points $P_1(1.988, 0.001)$, $P_2(4.5, 0.001)$, $Q_1(1.988, 2)$, $Q_2(4.5, 2)$ where the first coordinate represents the initial tumor radius $R(0)$ and the second represents the membrane rigidity parameter $\mathscr{S}^{-1}$. When $\mathcal{A} = 0.5$, linear theory predicts that a circular interface will evolve to its stationary radius $R_s \approx 3.326$, while the stability of a perturbed interface depends on $\mathscr{S}^{-1}$. As seen from Fig. \ref{fig:marginalSinv} [a], linear theory predicts that starting from the point $P_1(1.988, 0.001)$, $P_2(4.5,0.001)$, where membrane rigidity is low, the 3-mode perturbation will be unstable as the tumor grows or shrinks to the stationary radius $R_s$. We point out that the effective tumor radius will actually grow/shrink to the stationary radius $R_s(\mathcal{A}=0.5)\approx 3.326$ for a while and turn out to be larger than this predicted size due to morphological instability. However, starting from the point $Q_1(1.988, 2)$, $Q_2(4.5, 2)$ where the membrane rigidity is higher ($\mathscr{S}^{-1}=2$), the simulation shows tumor will grow or shrink to $R\approx3.326$ and is stable.

In Fig. \ref{fig:marginalSinv} [b] and [c], we show the nonlinear evolution of tumors starting from points $P_1$ and $Q_1$, respectively.  In each case, the initial tumor/host interface is a 3-mode perturbation of a circle given by $r=1.988+0.05\cos(3\alpha)$. Similar computations for $P_2$ and $Q_2$ are shown in Fig. \ref{fig:marginalSinv} [d] and [e] with initial shape $r=4.5+0.05\cos(3\alpha)$.  The viscosity ratio is varied from $\lambda=0.5, 1.5, 2.5$ as labeled in the plots. The nonlinear shape factors are plotted as functions of time, and the corresponding tumor morphology at the final time $t_f$ of each simulation are shown as insets. We also plot the linear solution for $\lambda=0.5$ (dashed line) in Fig. \ref{fig:marginalSinv} [b] and [d]. Note that the calculations for Fig. \ref{fig:marginalSinv} [b] and [d] stop because higher numerical resolution is needed to resolve high curvature regions.

When the membrane rigidity is small ($\mathscr{S}^{-1} = 0.001$), Fig. \ref{fig:marginalSinv} [b] show that the shape perturbation starting from the point $P_1$ grows rapidly, especially for small viscosity ratio, indicating an unstable growth.  This is consistent with the predictions of linear stability theory (Fig. \ref{fig:marginalSinv}[a]).  On the other hand, when the membrane rigidity is increased to $\mathscr{S}^{-1} = 2$, nonlinear simulations from point $Q_1$ converge to a final circular morphology for all three cases and  the tumor grows stably to its diffusion limited size  $R_s\approx3.326$, as shown in Fig. \ref{fig:marginalSinv}[c]. Simulations of points $P_2$ and $Q_2$ in Fig. \ref{fig:marginalSinv} [d] and [e] show similar unstable and stable shrinking behavior, respectively.

In Fig. \ref{fig:marginalSinv} [f], using the point $P_1$, we compare the shape perturbation between the non-weakening bending model (constant bending stiffness plotted using solid lines) and  curvature weakening bending model (plotted using dashed lines). Here we set $\lambda_c=1.25, C=0.95$ in Eq. \eqref{weakening}. For small viscosity ratio $\lambda=1.5$, the curvature weakening effect dramatically slows down the process of unstable growth and leads the interface morphology to branching patterns. However, as $\lambda$ increases, the weakening effect is reduced and we get the usual encapsulated morphology. These results highlight the level of complexity and sensitivity of interface dynamics in fluid due to inhomogeneous elasticity.  

\begin{figure}
\begin{minipage}{\linewidth}
    \centering
    \includegraphics[width=0.8\textwidth]{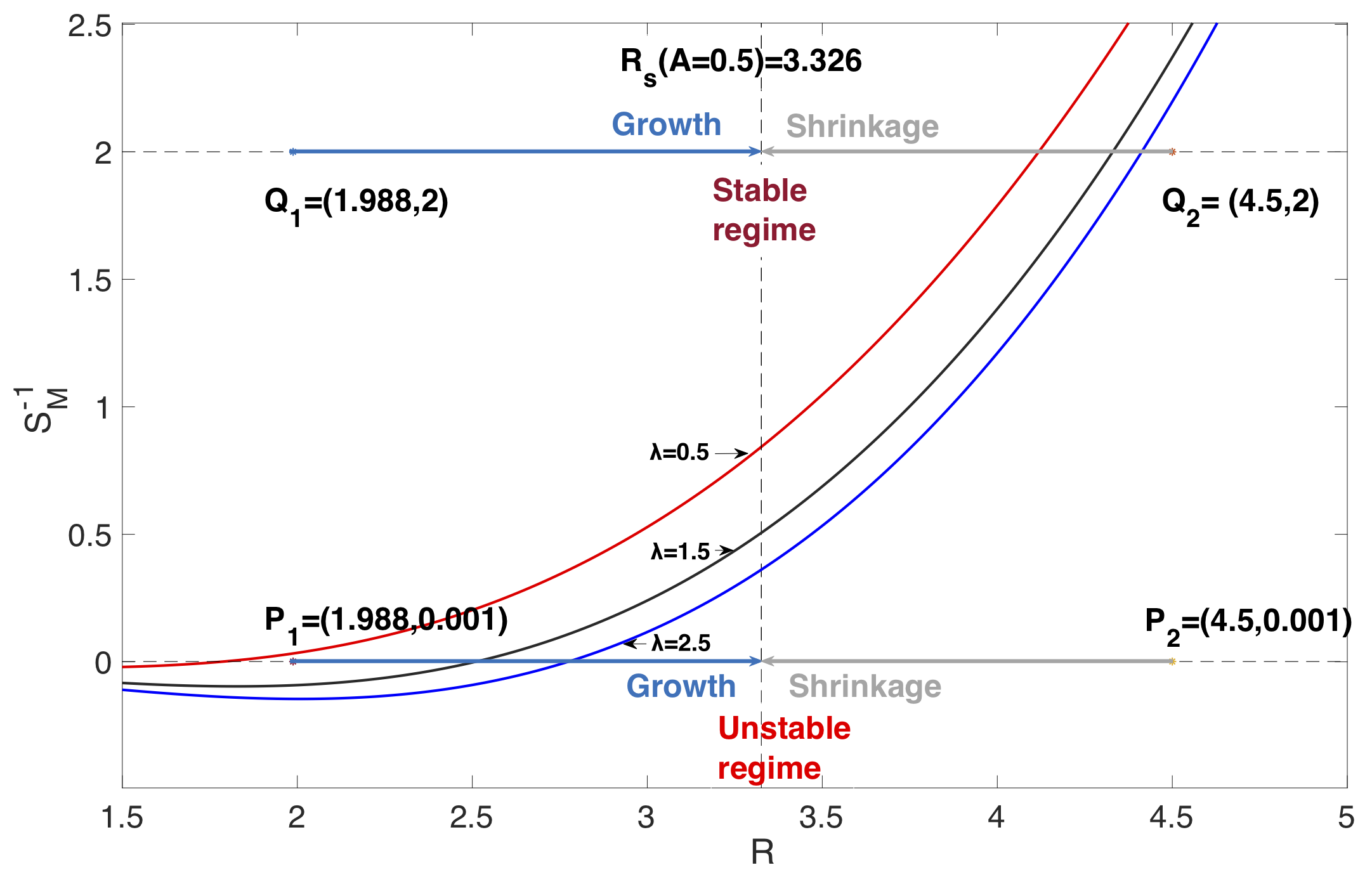}\textbf{[a]}
\end{minipage}
\begin{minipage}{\linewidth}
    \centering
    \includegraphics[width=0.8\textwidth]{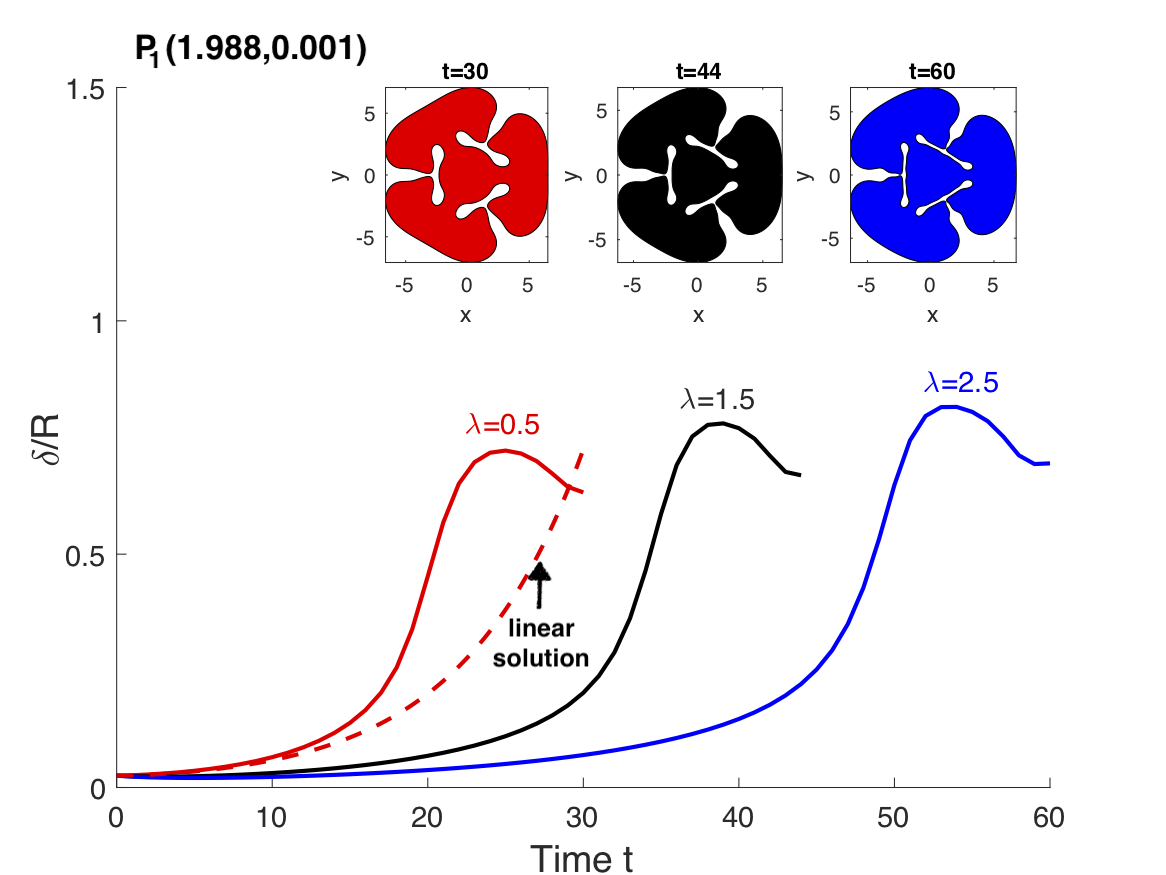}\textbf{[b]}
\end{minipage}
\begin{minipage}{\linewidth}
    \centering
    \includegraphics[width=0.8\textwidth]{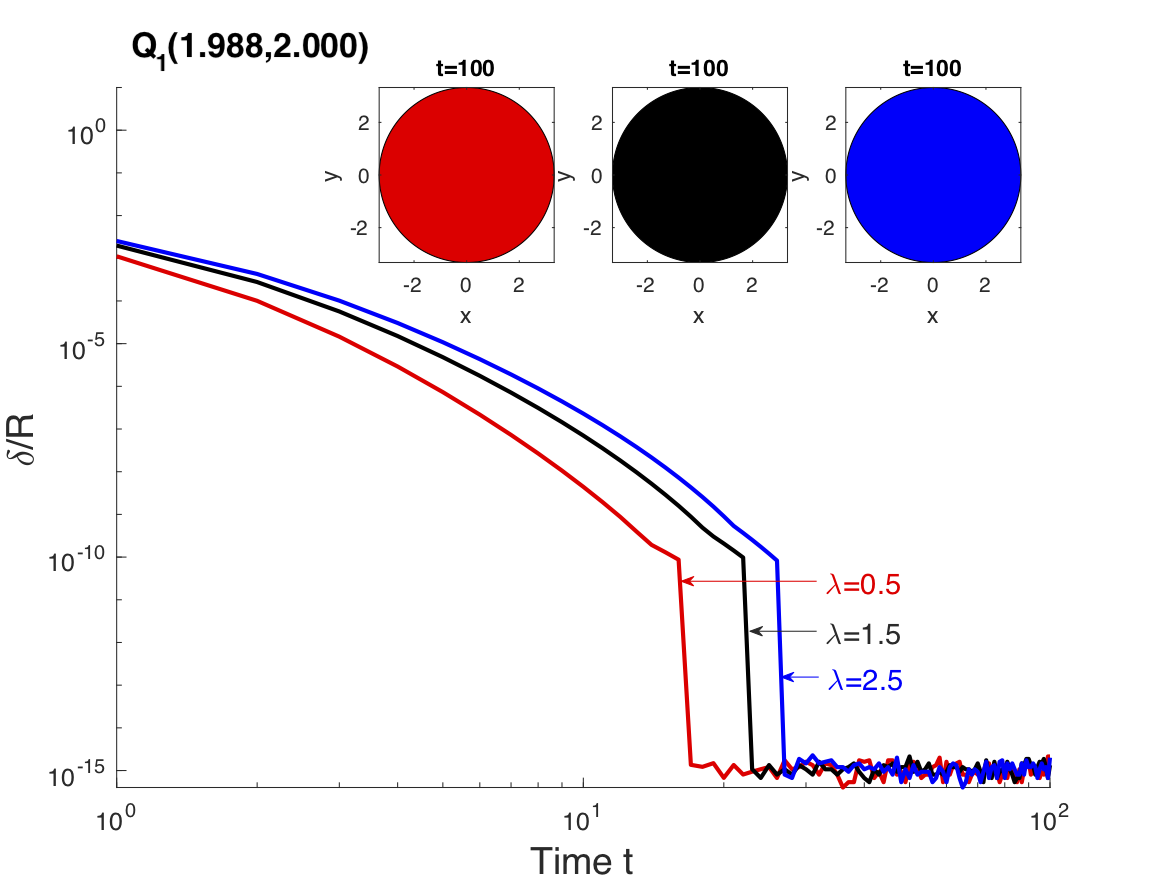}\textbf{[c]}
    \newpage
    \caption{[a] The marginally stable value of the membrane rigidity parameter $\mathscr{S}^{-1}$ as a function of unperturbed radius $R$ for different viscosity ratios (as labeled) with $A = 0.5$ and mode $l = 3$. The four points $P_1(1.988, 0.001), P_2(4.5, 0.001), Q_1(1.988, 2), Q_2(4.5, 2)$ indicate parameter values at which nonlinear simulations will be performed (see Fig. \ref{fig:marginalSinv} [b], [c], [d], and [e]). [b] Unstable growth corresponding to the point $P_1$ in Fig. \ref{fig:marginalSinv}[a]. [c] Stable growth corresponding to the point $Q_1$. The initial interface $r=1.988+0.05\cos(3\alpha)$ for $P_1$ and $Q_1$, and $r=4.5+0.05\cos(3\alpha)$ for $P_2$ and $Q_2$. We set $\mathscr{S}^{-1} = 0.001$  for $P_1$ and $P_2$, and $\mathscr{S}^{-1} = 2$  for $Q_1$ and $Q_2$. The viscosity ratios are $\lambda=0.5,1.5,2.5$ labeled with red, black and blue color respectively.}
    \label{fig:marginalSinv}
    \end{minipage}
\end{figure}

\setcounter{figure}{2}  

\begin{figure}
\begin{minipage}{\linewidth}
    \centering
    \includegraphics[width=0.8\textwidth]{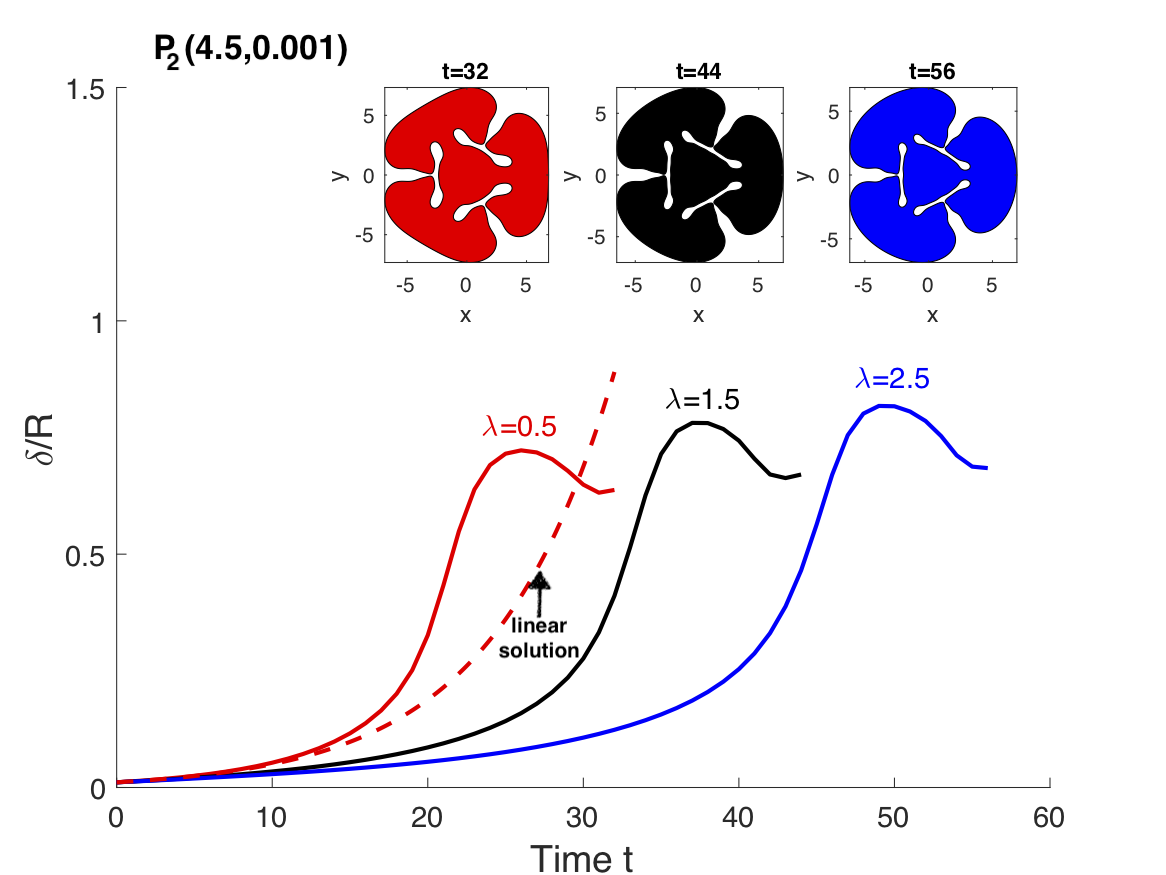}\textbf{[d]}
\end{minipage}
\begin{minipage}{\linewidth}
    \centering
    \includegraphics[width=0.8\textwidth]{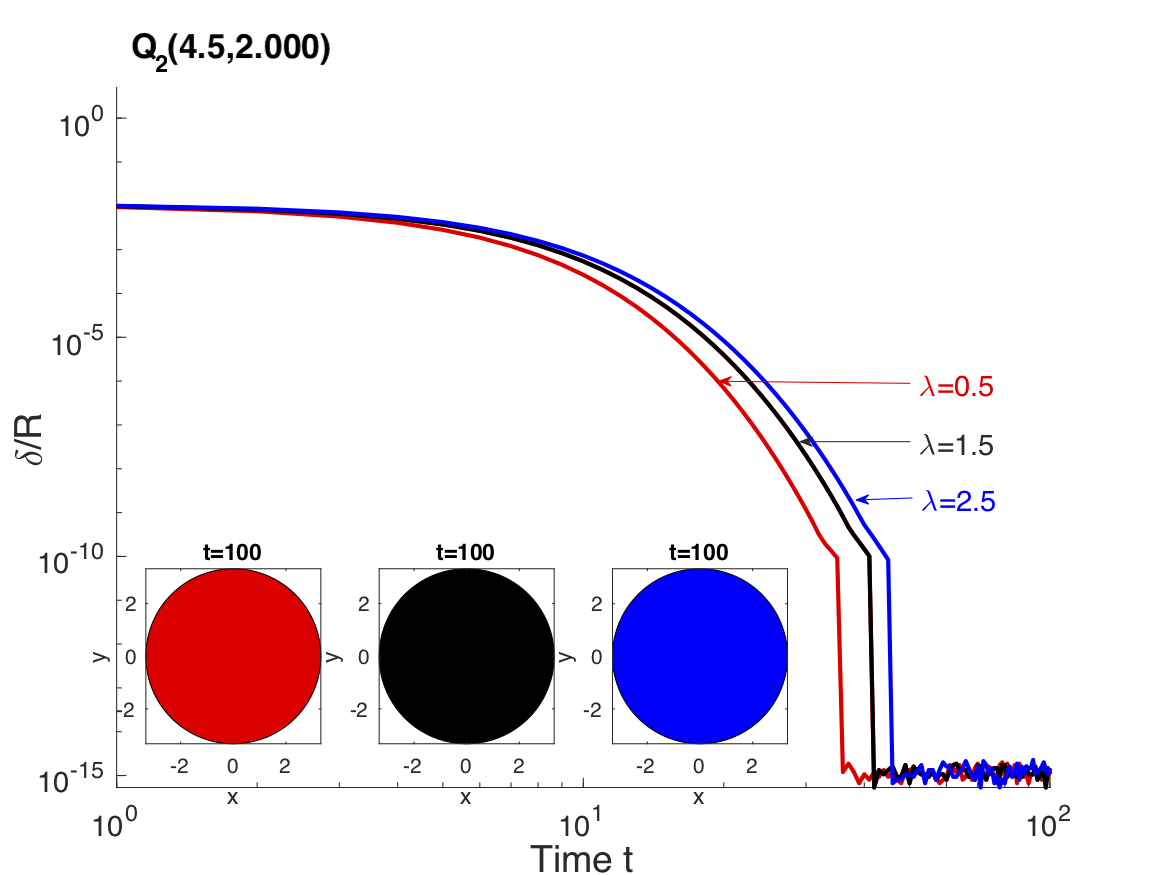}\textbf{[e]}
\end{minipage}
\begin{minipage}{\linewidth}
    \centering
    \includegraphics[width=0.8\textwidth]{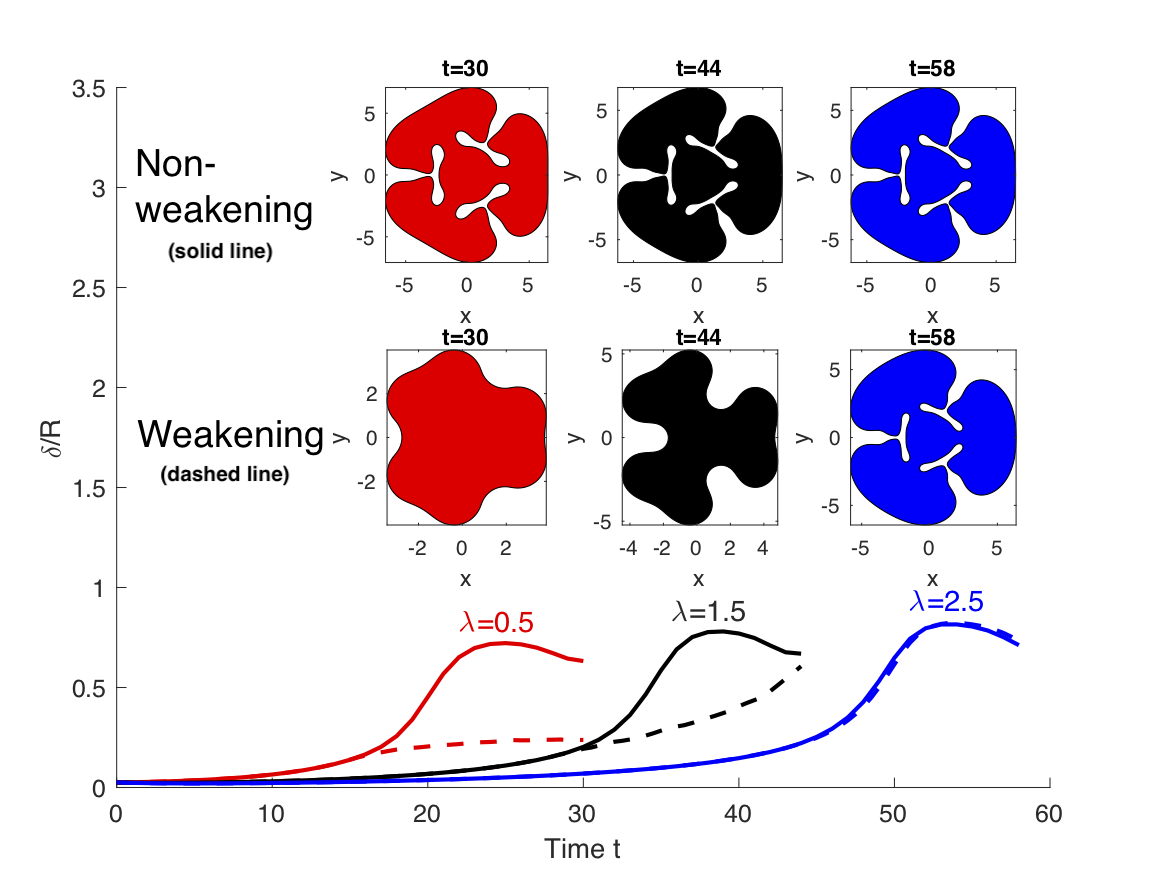}\textbf{[f]}
    \caption{[d] Unstable growth corresponding to the point $P_2$. [e] Stable growth corresponding to the point $Q_2$. [f] Comparison between the uniform bending and curvature weakening bending.}
    \label{fig:marginalSinv}
    \end{minipage}
\end{figure}

\subsection{Growth in nonlinear regime ($\mathcal{A}=0.7$)}
\begin{figure}
\begin{minipage}{\linewidth}
    \centering
    \includegraphics[width=0.8\textwidth]{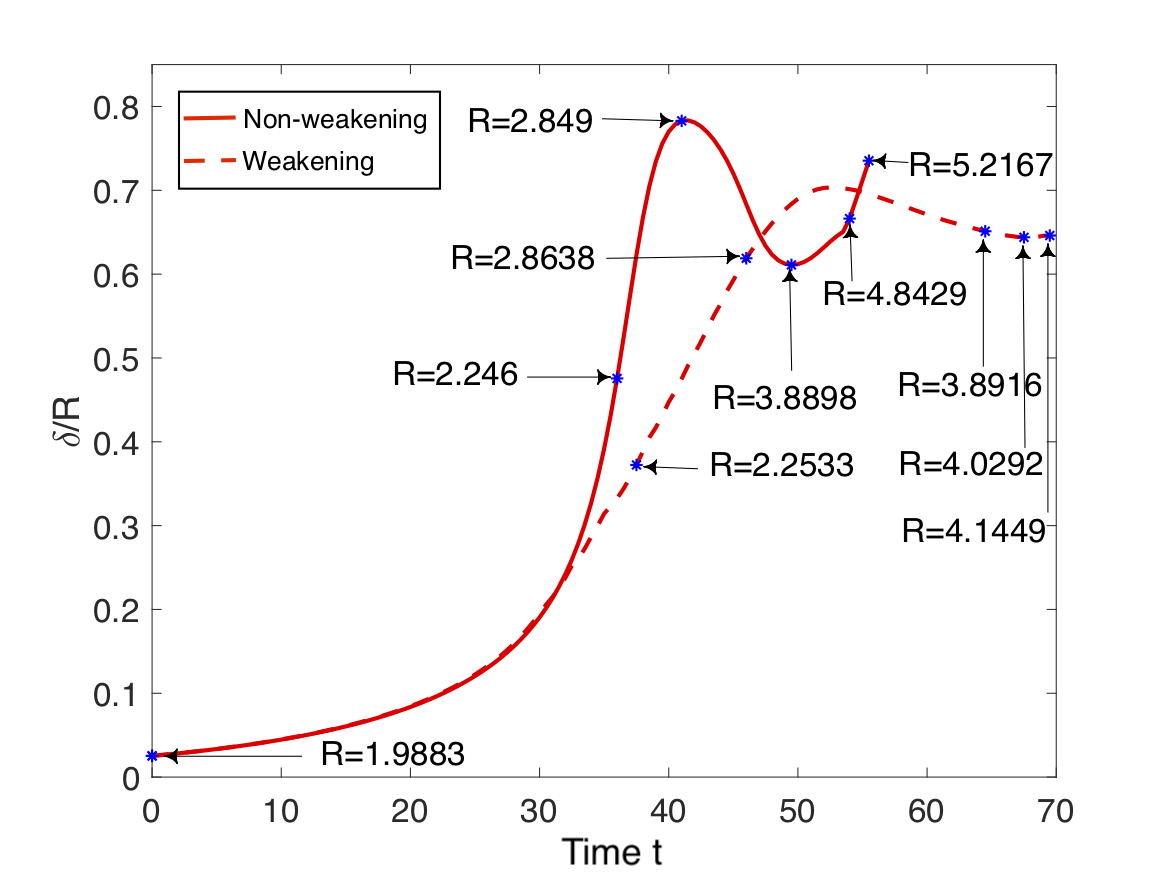}\textbf{[a]}
\end{minipage}
\begin{minipage}{\linewidth}
    \centering
    \includegraphics[width=0.8\textwidth]{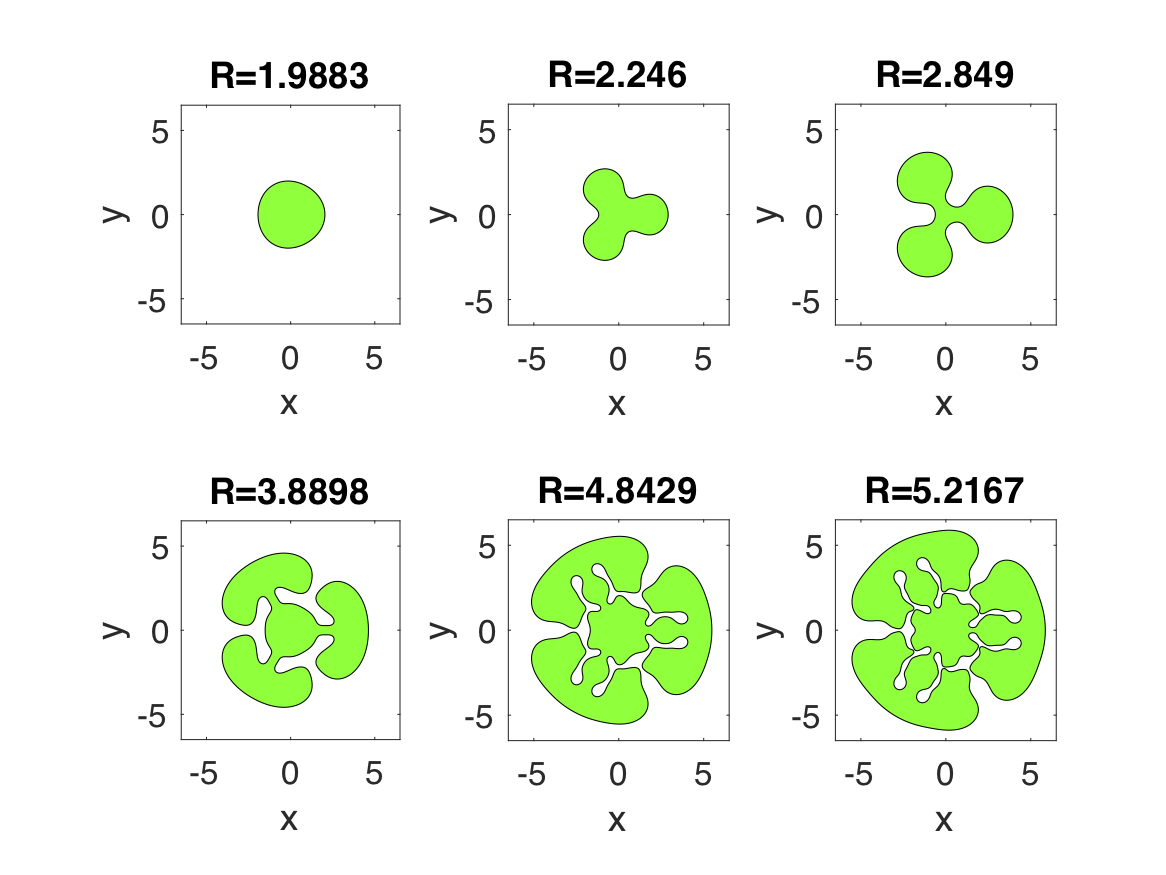}\textbf{[b]}
\end{minipage}
\begin{minipage}{\linewidth}
    \centering
    \includegraphics[width=0.8\textwidth]{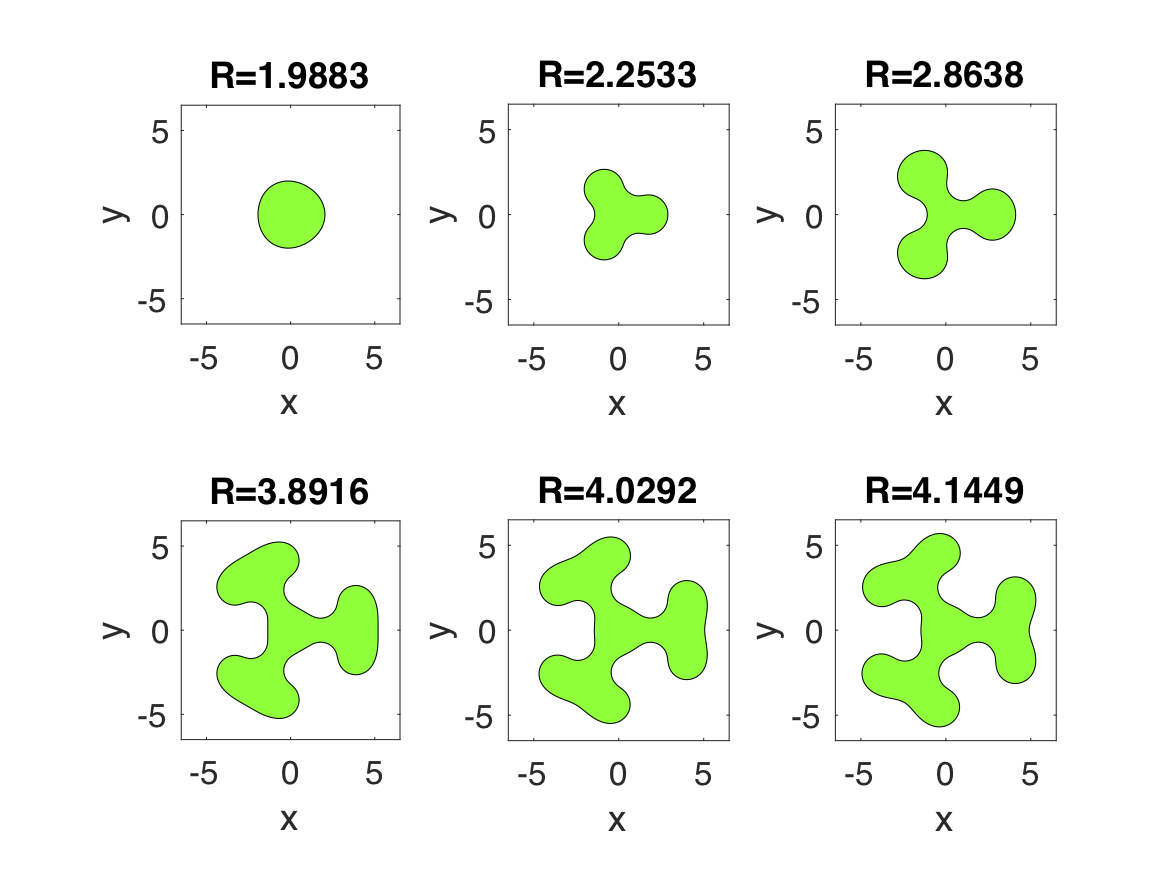}\textbf{[c]}
    \caption{Shape factor evolution[a] of the growth of Stokes-flow model with non-weakening[b]/weakening[c] bending energy beyond the linear prediction with a simple tumor initial shape; the parameters are set as: mesh points $N=2048$, time step $dt=0.01$, bending rigidity $\mathscr{S}^{-1}=0.001$, viscosity ratio $\lambda=1.5$, apoptosis $\mathscr{A}=0.7$, initial shape: $r=1.988+0.05\cos3\alpha$.}
    \label{fig:nonlinear regime simple}
\end{minipage}
\end{figure}

 We next increase the apoptosis rate $\mathcal{A}$ from 0.5 to 0.7 and focus only on the viscosity ratio $\lambda=1.5$ case. Beyond the parameter regime of linear prediction, in Fig. \ref{fig:nonlinear regime simple} [a] the interface evolves far away from the steady radius $R_s(\mathcal{A}=0.7)=1.988$. Notice that in this case based on its initial 3-mode perturbation, i.e. $r=1.988+0.05\cos(3\alpha)$, the interface develops further splitting and fingering patterns due to the curvature weakening effect as shown in Fig. \ref{fig:nonlinear regime simple} [c].  The evolving morphology changes from a compact shape with inward splitting shown in Fig.  \ref{fig:nonlinear regime simple} [b] to a fingering pattern with a tendency of outward splitting as shown in Fig. \ref{fig:nonlinear regime simple} [c]. Here we use curvature weakening parameters $\lambda_c=1.25,\ C=0.95$.

\begin{figure}
\begin{minipage}{\linewidth}
    \centering
    \includegraphics[width=\textwidth]{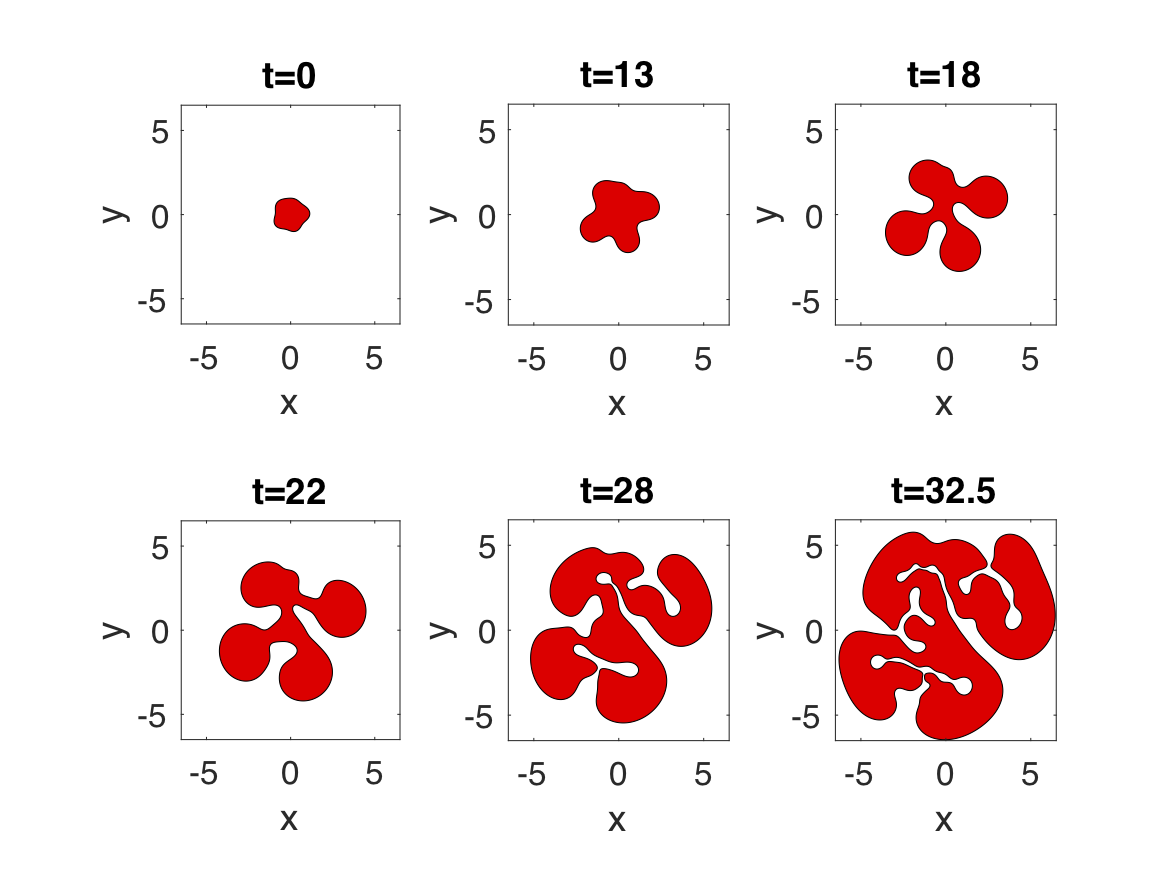}\textbf{[a]}
\end{minipage}
\begin{minipage}{\linewidth}
    \centering
    \includegraphics[width=\textwidth]{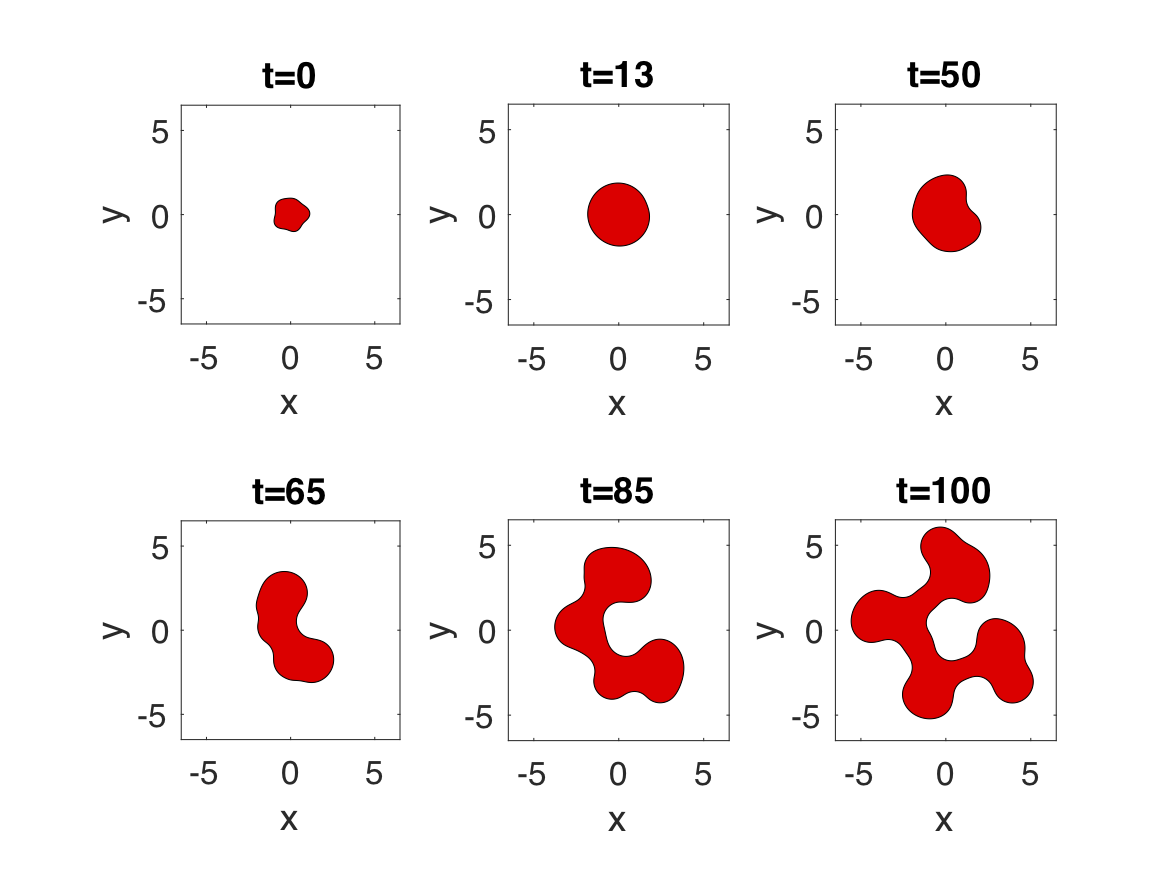}\textbf{[b]}
    \caption{Growth of Stokes-flow model with non-weakening[a]/weakening($\lambda_c=1.25,C=0.95$)[b] bending energy beyond the linear prediction with a complex initial shape; the parameters are set as: mesh points $N=2048$, time step $dt=0.01$, bending rigidity $\mathscr{S}^{-1}=0.001$, viscosity ratio $\lambda=1.5$, apoptosis $\mathscr{A}=0.7$, initial shape: $r=1+\frac{0.05}{1.988} \cos (2 \alpha)+\frac{0.1}{1.988} \cos (3 \alpha)+\frac{0.08}{1.988} \sin (4 \alpha) +\frac{0.12}{1.988} \cos (5 \alpha)$.}
    \label{fig:nonlinear regime complex}
\end{minipage}
\end{figure}

We also compute the evolution of a complex initial shape in Fig. \ref{fig:nonlinear regime complex}, where $\displaystyle r=1+\frac{0.05}{1.988} \cos (2 \alpha)+\frac{0.1}{1.988} \cos (3 \alpha)+\frac{0.08}{1.988} \sin (4 \alpha) +\frac{0.12}{1.988} \cos (5 \alpha).$
In Fig. \ref{fig:nonlinear regime complex} [a], tumor evolution  with complex initial shape grows unstably with appearance of long and slim zigzags inside the tumor; while for the one with curvature weakening effect ($\lambda_c=1.25,\ C=0.95$) in Fig. \ref{fig:nonlinear regime complex} [b] such pattern disappears, and it takes much longer time to reach the size in Fig. \ref{fig:nonlinear regime complex} [a].

\subsection{Self-similar patterns (Time-varying $\mathcal{A}$)}
\begin{figure}
\begin{minipage}{\linewidth}
    \centering
    \includegraphics[width=\textwidth]{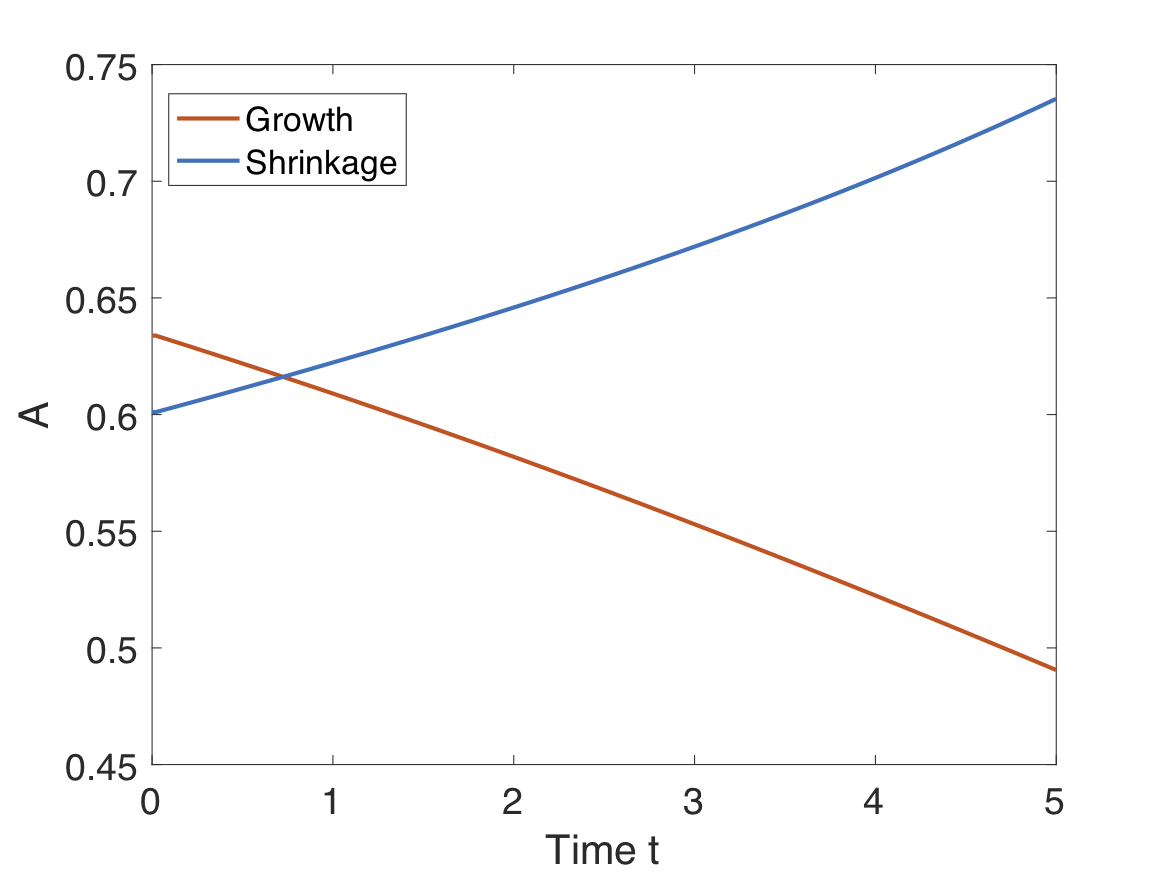}\textbf{[a]}
    \end{minipage}
\begin{minipage}{\linewidth}
    \centering
    \includegraphics[width=\textwidth]{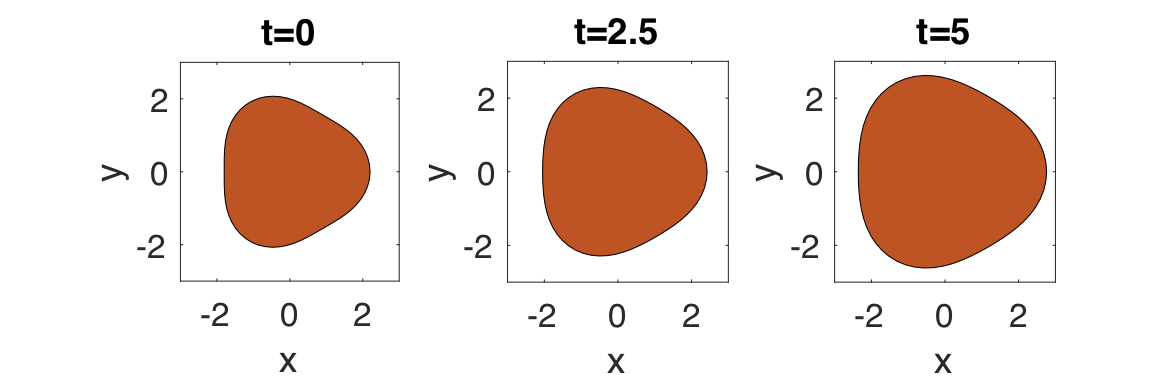}\textbf{[b]}
    \end{minipage}
\begin{minipage}{\linewidth}
    \centering
    \includegraphics[width=\textwidth]{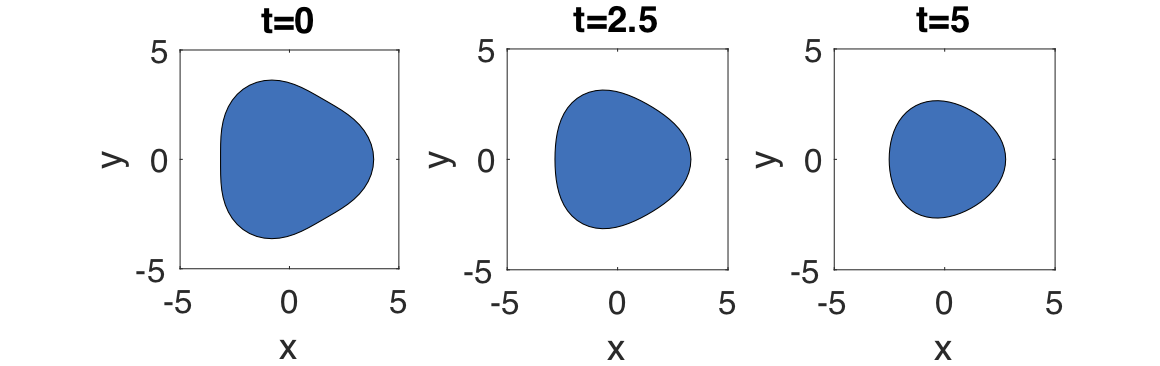}\textbf{[c]}
    \caption{[a] The time evolution of the relative rate of cell apoptosis to mitosis $\mathcal{A}$ to keep the self-similar shape of tumor in [b] and [c]. [b]/[c] Self-similar growth[b]/shrinkage[c] of Stokes-flow model with non-weakening bending energy by the prediction of morphological stability in Eq. \eqref{shape factor}. Here we set viscosity ratio $\lambda=0.5$ for [b] and $7.5$ for [c]. Apoptosis $\mathcal{A}$ is time-varying such that the shape factor $\frac{\delta}{R}=0$ in linear regime.  Initial shape: $r=2+0.2\cos3\alpha$ for [b] and $r=3.5+0.35\cos3\alpha$ for [c].}
    \label{fig:selfsimilar}
    \end{minipage}
\end{figure}

 In Fig. \ref{fig:selfsimilar}, we demonstrate the existence of linear self-similar growth/shrinkage of a tumor with a 3-fold perturbation. Here we choose $\mathcal{A}$ such that the evolution of the shape factor $\displaystyle\frac{d(\frac{\delta}{R})}{dt}=0$, following Eq. \eqref{shape factor}. That is  $\frac{\delta}{R}(t)=\frac{\delta}{R}(0)$. To get this shape preserving evolution, the relative rate of cell apoptosis to mitosis $\mathcal{A}$ must be time (or size) dependent. As plotted in Fig. \ref{fig:selfsimilar} [a], $\mathcal{A}$ is decreasing/increasing function for growth/shrinkage, respectively. Snapshots of a self-similar sequence are plotted in Fig. \ref{fig:selfsimilar} [b] and [c].  Although preliminary, the self-similar idea helps shed light on the strategy for morphological control, as a time dependent apoptosis might be enforced by a well-designed chemo- or radiotherapy.
 
\section{Conclusion}
\label{conclusion}
In this paper, we performed  nonlinear simulations of a 2D, non-circular tumor with isotropic or  curvature weakened bending rigidity growing in a host tissue. The interior tumor and exterior host were modeled by the Stokes flow, and the tumor-host interface was modeled by an elastic membrane governed by the Helfrich bending energy. Using boundary integral formulations of the Stokes flow and the nutrient field, we developed a spectrally accurate sharp interface method. We then investigated the nonlinear dynamics of the tumor-host interface.

\begin{figure}
\begin{minipage}{\linewidth}
    \centering
    \includegraphics[width=\textwidth]{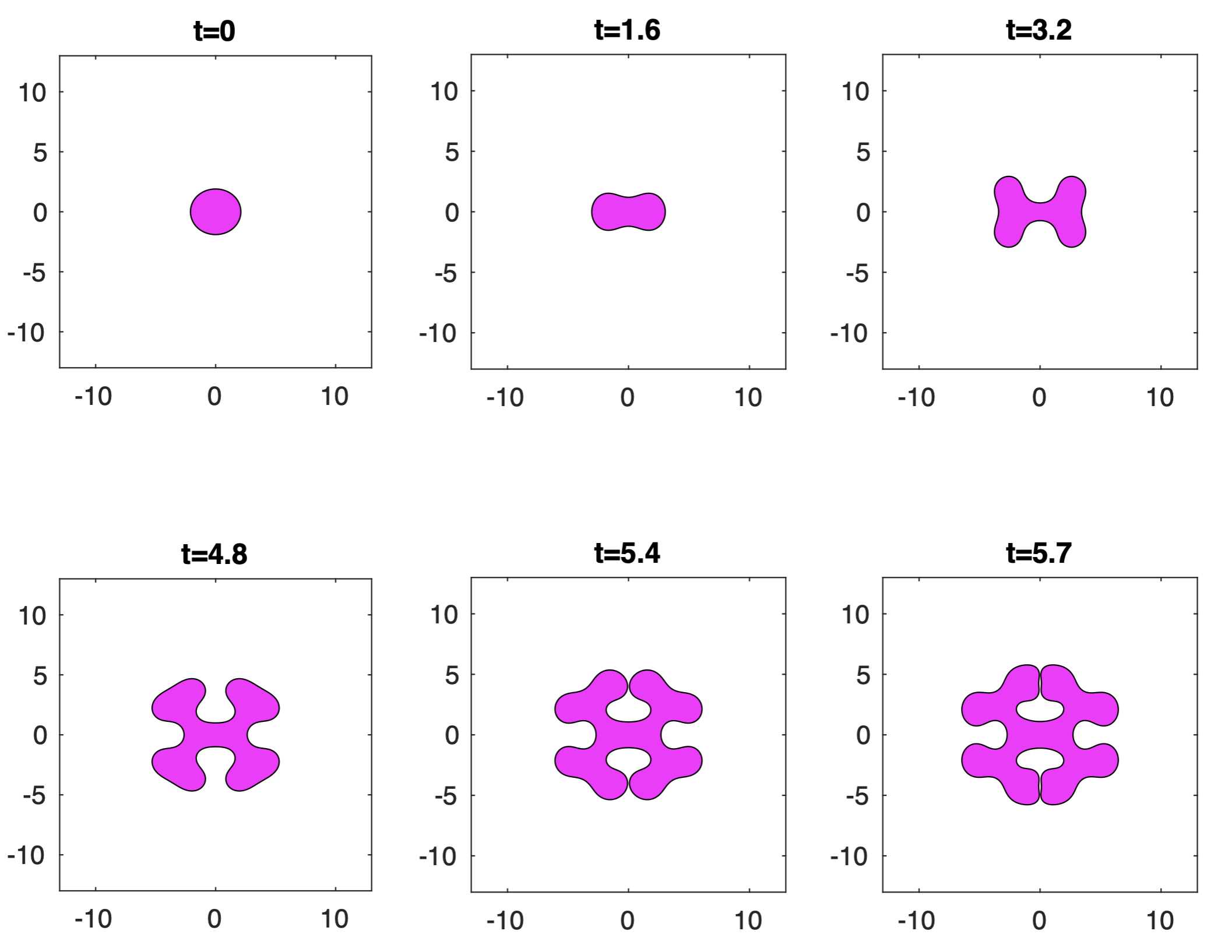}\textbf{[a]}
    \end{minipage}
    \begin{minipage}{\linewidth}
    \centering
    \includegraphics[width=\textwidth]{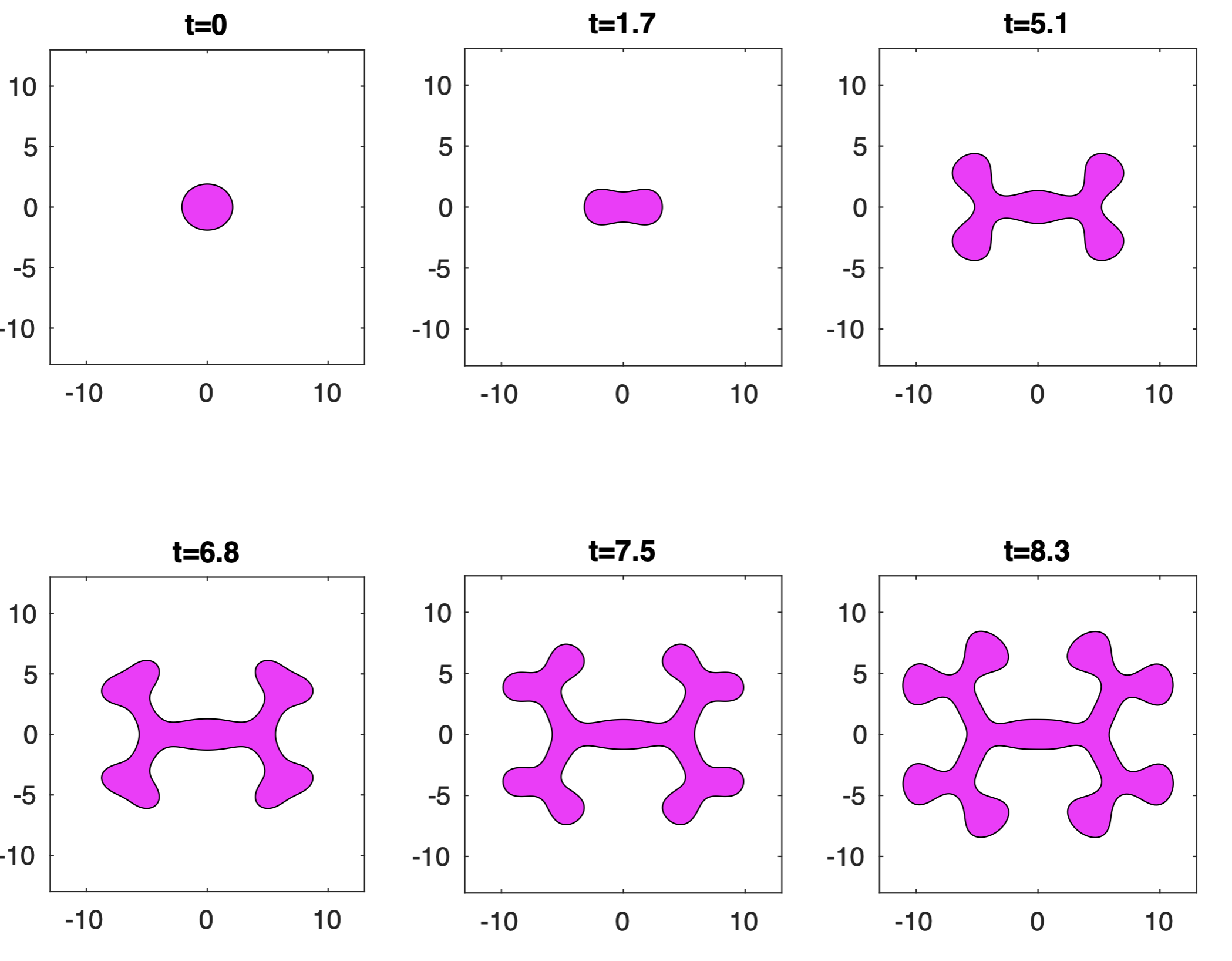}\textbf{[b]}
    \end{minipage}
    \caption{Growth with Darcy-flow model with [a] isotropic bending energy and [b] curvature weakening bending (stiffness fraction $C=1$, characteristic length $\lambda_c=1$); the parameters are set as: mesh points $N=2048$, time step $dt=0.01$, bending rigidity $\mathscr{S}^{-1}=0.15$, apoptosis $\mathscr{A}=0.7$, initial shape: $r=2+0.05 \cos (2 \alpha)$.}
    \label{fig:darcyS0.15}
\end{figure}

The linear stability analysis suggests that an increase in bending rigidity contributes to an increase in morphological stability for an isotropic bending rigidity. Nonlinear simulation confirms this and moreover, curvature weakening bending helps improve the stability by slowing down the growth of shape perturbations and promotes branching or tip-splitting fingering patterns rather than encapsulated morphologies for small viscosity ratio, In fact, not only for the Stokes model, our recent preliminary results using the Darcy's model suggest more pronounced fingering patterns if a curvature weakening bending is implemented, as shown in Fig. \ref{fig:darcyS0.15}. We can see the self-branching morphology is enhanced with curvature bending energy as reported in \cite{zhao2016nonlinear} for a Hele-Shaw interface. It is also observed that an increase in the apoptosis leads to an overall increase in shape instabilities. 

In experiments, thermal and mechanical stresses have been found to be important in regulating cell fates and motility, proliferation and apoptosis rates. In future work, we will consider these effects. We will also consider adding a stochastic component to the current model. Although our 2D results are expected to hold qualitatively in three dimensions as suggested by the linear stability analysis (at least for Darcy's model), we would like to perform full 3D simulations to confirm this.

 {We have chosen the initial shape as a perturbed circle since that that in vitro tumor grows nearly spherical at early times.  It is reasonable to assume that, in vivo, tumor at its initial stage of avascular growth is nearly spherical. As the tumor continues to grow, morphological instabilities come in and the tumor starts to develop protruding shapes.  In this work we focus on the morphological instability of the interface and our numerical scheme can handle even complex tumor morphologies as indicated in Fig. 4[a] and Fig. 5[a] with spectral accuracy in space. For studies involving simulation beyond topological changes, we plan to use phase-field or level-set formulation, which is our future work.}

\begin{acknowledgements}
S. L. acknowledges the support from the National Science Foundation, Division of Mathematical Sciences grant DMS-1720420. S. L. was also partially supported by grant ECCS-1307625. M. L. acknowledges the F. R. “Buck” McMorris Summer Research support from the College of Science, IIT. C. L. is partially supported by the National Science Foundation, Division of Mathematical Sciences grant DMS-1759536.

\end{acknowledgements}

\bibliographystyle{spr-chicago}      
\bibliography{template}   

\begin{thebibliography}{51}
\ifx \bisbn   \undefined \def \bisbn  #1{ISBN #1}\fi
\ifx \binits  \undefined \def \binits#1{#1} \fi
\ifx \bauthor  \undefined \def \bauthor#1{#1} \fi
\ifx \bjtitle  \undefined \def \bjtitle#1{\textit{#1}}\fi
\ifx \batitle  \undefined \def \batitle#1{#1} \fi
\ifx \bctitle  \undefined \def \bctitle#1{#1} \fi
\ifx \bvolume  \undefined \def \bvolume#1{#1}\fi
\ifx \byear  \undefined \def \byear#1{#1} \fi
\ifx \bissue  \undefined \def \bissue#1{#1} \fi
\ifx \bfpage  \undefined \def \bfpage#1{#1} \fi
\ifx \blpage  \undefined \def \blpage #1{#1} \fi
\ifx \burl  \undefined \def \burl#1{#1} \fi
\ifx \doiurl  \undefined \def \doiurl#1{#1} \fi
\ifx \betal  \undefined \def \betal{et al.} \fi
\ifx \binstitute  \undefined \def \binstitute#1{#1} \fi
\ifx \beditor  \undefined \def \beditor#1{#1} \fi
\ifx \bpublisher  \undefined \def \bpublisher#1{#1} \fi
\ifx \bbtitle  \undefined \def \bbtitle#1{\textit{#1}} \fi
\ifx \bedition  \undefined \def \bedition#1{#1} \fi
\ifx \bseriesno  \undefined \def \bseriesno#1{#1} \fi
\ifx \blocation  \undefined \def \blocation#1{#1} \fi
\ifx \bsertitle  \undefined \def \bsertitle#1{\textit{#1}} \fi
\ifx \bsnm \undefined \def \bsnm#1{#1} \fi
\ifx \bsuffix \undefined \def \bsuffix#1{#1} \fi
\ifx \bparticle \undefined \def \bparticle#1{#1} \fi
\ifx \barticle \undefined \def \barticle#1{#1} \fi
\ifx \botherref \undefined \def \botherref #1{#1} \fi
\ifx \url \undefined \def \url#1{#1} \fi
\ifx \bchapter \undefined \def \bchapter#1{#1} \fi
\ifx \bbook \undefined \def \bbook#1{#1} \fi
\ifx \bcomment \undefined \def \bcomment#1{#1} \fi
\ifx \oauthor \undefined \def \oauthor#1{#1} \fi
\ifx \citeauthoryear \undefined \def \citeauthoryear#1{#1} \fi
\ifx \texttildelow  \undefined \def \texttildelow{\symbol{126}} \fi
\def \endbibitem {}
\ifx \bconflocation  \undefined \def \bconflocation#1{#1} \fi

\bibitem[\protect\citeauthoryear{Adkins and Zhou}{2017}]{Zhou20172}
\begin{barticle}
\bauthor{\bsnm{Adkins}, \binits{Melissa~R}}, and \bauthor{\binits{YC}
  \bsnm{Zhou}}.
\byear{2017}.
\batitle{Geodesic curvature driven surface microdomain formation}.
\bjtitle{Journal of computational physics}
\bvolume{345}: \bfpage{260}--\blpage{274}.
doi:\doiurl{10.1016/j.jcp.2017.05.029}.
\end{barticle}
\endbibitem

\bibitem[\protect\citeauthoryear{Baker and Shelley}{1990}]{baker1990connection}
\begin{barticle}
\bauthor{\bsnm{Baker}, \binits{GR}}, and \bauthor{\binits{MJ} \bsnm{Shelley}}.
\byear{1990}.
\batitle{On the connection between thin vortex layers and vortex sheets}.
\bjtitle{Journal of Fluid Mechanics}
\bvolume{215}: \bfpage{161}--\blpage{194}.
\end{barticle}
\endbibitem

\bibitem[\protect\citeauthoryear{Bellomo and
  de~Angelis}{2008}]{bellomo2008selected}
\begin{bbook}
\bauthor{\bsnm{Bellomo}, \binits{Nicola}}, and \bauthor{\binits{Elena}
  \bparticle{de~}\bsnm{Angelis}}.
\byear{2008}.
\bbtitle{Selected topics in cancer modeling: genesis, evolution, immune
  competition, and therapy}.
\bpublisher{Springer}.
\end{bbook}
\endbibitem

\bibitem[\protect\citeauthoryear{Byrne and Chaplain}{1996}]{byrne1996growth}
\begin{barticle}
\bauthor{\bsnm{Byrne}, \binits{H~M\_}}, and \bauthor{\binits{MAJ}
  \bsnm{Chaplain}}.
\byear{1996}.
\batitle{Growth of necrotic tumors in the presence and absence of inhibitors}.
\bjtitle{Mathematical biosciences}
\bvolume{135} (\bissue{2}): \bfpage{187}--\blpage{216}.
\end{barticle}
\endbibitem

\bibitem[\protect\citeauthoryear{Byrne}{2010}]{byrne2010dissecting}
\begin{barticle}
\bauthor{\bsnm{Byrne}, \binits{Helen~M}}.
\byear{2010}.
\batitle{Dissecting cancer through mathematics: from the cell to the animal
  model}.
\bjtitle{Nature Reviews Cancer}
\bvolume{10} (\bissue{3}): \bfpage{221}.
\end{barticle}
\endbibitem

\bibitem[\protect\citeauthoryear{Byrne and Chaplain}{1997}]{byrne1997free}
\begin{barticle}
\bauthor{\bsnm{Byrne}, \binits{HM}}, and \bauthor{\binits{Mark~AJ}
  \bsnm{Chaplain}}.
\byear{1997}.
\batitle{Free boundary value problems associated with the growth and
  development of multicellular spheroids}.
\bjtitle{European Journal of Applied Mathematics}
\bvolume{8} (\bissue{6}): \bfpage{639}--\blpage{658}.
\end{barticle}
\endbibitem

\bibitem[\protect\citeauthoryear{Chen et~al.}{2014}]{chen2014stable}
\begin{barticle}
\bauthor{\bsnm{Chen}, \binits{Ying}}, \bauthor{\binits{Steven~M} \bsnm{Wise}},
  \bauthor{\binits{Vivek~B} \bsnm{Shenoy}}, and \bauthor{\binits{John~S}
  \bsnm{Lowengrub}}.
\byear{2014}.
\batitle{A stable scheme for a nonlinear, multiphase tumor growth model with an
  elastic membrane}.
\bjtitle{International journal for numerical methods in biomedical engineering}
\bvolume{30} (\bissue{7}): \bfpage{726}--\blpage{754}.
\end{barticle}
\endbibitem

\bibitem[\protect\citeauthoryear{Cristini and
  Lowengrub}{2010}]{cristini2010multiscale}
\begin{bbook}
\bauthor{\bsnm{Cristini}, \binits{Vittorio}}, and \bauthor{\binits{John}
  \bsnm{Lowengrub}}.
\byear{2010}.
\bbtitle{Multiscale modeling of cancer: an integrated experimental and
  mathematical modeling approach}.
\bpublisher{Cambridge University Press}.
\end{bbook}
\endbibitem

\bibitem[\protect\citeauthoryear{Cristini
  et~al.}{2017}]{cristini2017introduction}
\begin{bbook}
\bauthor{\bsnm{Cristini}, \binits{Vittorio}}, \bauthor{\binits{Eugene}
  \bsnm{Koay}}, and \bauthor{\binits{Zhihui} \bsnm{Wang}}.
\byear{2017}.
\bbtitle{An introduction to physical oncology: How mechanistic mathematical
  modeling can improve cancer therapy outcomes}.
\bpublisher{CRC Press}.
\end{bbook}
\endbibitem

\bibitem[\protect\citeauthoryear{Cristini et~al.}{2003}]{cristini2003nonlinear}
\begin{barticle}
\bauthor{\bsnm{Cristini}, \binits{Vittorio}}, \bauthor{\binits{John}
  \bsnm{Lowengrub}}, and \bauthor{\binits{Qing} \bsnm{Nie}}.
\byear{2003}.
\batitle{Nonlinear simulation of tumor growth}.
\bjtitle{Journal of mathematical biology}
\bvolume{46} (\bissue{3}): \bfpage{191}--\blpage{224}.
\end{barticle}
\endbibitem

\bibitem[\protect\citeauthoryear{Cristini
  et~al.}{2005}]{cristini2005morphologic}
\begin{barticle}
\bauthor{\bsnm{Cristini}, \binits{Vittorio}}, \bauthor{\binits{Hermann~B}
  \bsnm{Frieboes}}, \bauthor{\binits{Robert} \bsnm{Gatenby}},
  \bauthor{\binits{Sergio} \bsnm{Caserta}}, \bauthor{\binits{Mauro}
  \bsnm{Ferrari}}, and \bauthor{\binits{John} \bsnm{Sinek}}.
\byear{2005}.
\batitle{Morphologic instability and cancer invasion}.
\bjtitle{Clinical Cancer Research}
\bvolume{11} (\bissue{19}): \bfpage{6772}--\blpage{6779}.
\end{barticle}
\endbibitem

\bibitem[\protect\citeauthoryear{Dai and Promislow}{2013}]{Dai2013}
\begin{barticle}
\bauthor{\bsnm{Dai}, \binits{Shibin}}, and \bauthor{\binits{Keith}
  \bsnm{Promislow}}.
\byear{2013}.
\batitle{Geometric evolution of bilayers under the functionalized
  cahn--hilliard equation}.
\bjtitle{Proc. R. Soc. A.}
\bvolume{469}: \bfpage{20120505}.
\end{barticle}
\endbibitem

\bibitem[\protect\citeauthoryear{Du et~al.}{2004}]{DU2004450}
\begin{barticle}
\bauthor{\bsnm{Du}, \binits{Qiang}}, \bauthor{\binits{Chun} \bsnm{Liu}}, and
  \bauthor{\binits{Xiaoqiang} \bsnm{Wang}}.
\byear{2004}.
\batitle{A phase field approach in the numerical study of the elastic bending
  energy for vesicle membranes}.
\bjtitle{Journal of Computational Physics}
\bvolume{198} (\bissue{2}): \bfpage{450}--\blpage{468}.
\end{barticle}
\endbibitem

\bibitem[\protect\citeauthoryear{Du et~al.}{2005}]{Du2005}
\begin{barticle}
\bauthor{\bsnm{Du}, \binits{Qiang}}, \bauthor{\binits{Chun} \bsnm{Liu}},
  \bauthor{\binits{Rolf} \bsnm{Ryham}}, and \bauthor{\binits{Xiaoqiang}
  \bsnm{Wang}}.
\byear{2005}.
\batitle{A phase field formulation of the willmore problem}.
\bjtitle{Nonlinearity}
\bvolume{18}: \bfpage{1249}--\blpage{1267}.
\end{barticle}
\endbibitem

\bibitem[\protect\citeauthoryear{Frieboes
  et~al.}{2006}]{frieboes2006integrated}
\begin{barticle}
\bauthor{\bsnm{Frieboes}, \binits{Hermann~B}}, \bauthor{\binits{Xiaoming}
  \bsnm{Zheng}}, \bauthor{\binits{Chung-Ho} \bsnm{Sun}},
  \bauthor{\binits{Bruce} \bsnm{Tromberg}}, \bauthor{\binits{Robert}
  \bsnm{Gatenby}}, and \bauthor{\binits{Vittorio} \bsnm{Cristini}}.
\byear{2006}.
\batitle{An integrated computational/experimental model of tumor invasion}.
\bjtitle{Cancer research}
\bvolume{66} (\bissue{3}): \bfpage{1597}--\blpage{1604}.
\end{barticle}
\endbibitem

\bibitem[\protect\citeauthoryear{Friedman and Hu}{2007a}]{Hu20071}
\begin{barticle}
\bauthor{\bsnm{Friedman}, \binits{A.}}, and \bauthor{\binits{B.} \bsnm{Hu}}.
\byear{2007}a.
\batitle{Bifurcation for a free boundary problem modeling tumor growth by
  stokes equation}.
\bjtitle{SIAM Journal on Mathematical Analysis}
\bvolume{39} (\bissue{1}): \bfpage{174}--\blpage{194}.
doi:\doiurl{10.1137/060656292}.
\end{barticle}
\endbibitem

\bibitem[\protect\citeauthoryear{Friedman and Hu}{2007b}]{Hu20072}
\begin{barticle}
\bauthor{\bsnm{Friedman}, \binits{Avner}}, and \bauthor{\binits{Bei}
  \bsnm{Hu}}.
\byear{2007}b.
\batitle{Bifurcation from stability to instability for a free boundary problem
  modeling tumor growth by stokes equation}.
\bjtitle{Journal of mathematical analysis and applications}
\bvolume{327} (\bissue{1}): \bfpage{643}--\blpage{664}.
\end{barticle}
\endbibitem

\bibitem[\protect\citeauthoryear{Gavish et~al.}{2012}]{Promislow2012}
\begin{barticle}
\bauthor{\bsnm{Gavish}, \binits{Nir}}, \bauthor{\binits{Jaylan} \bsnm{Jones}},
  \bauthor{\binits{Zhengfu} \bsnm{Xu}}, \bauthor{\binits{Andrew}
  \bsnm{Christlieb}}, and \bauthor{\binits{Keith} \bsnm{Promislow}}.
\byear{2012}.
\batitle{Variational models of network formation and ion transport:
  Applications to perfluorosulfonate ionomer membranes}.
\bjtitle{Polymers}
\bvolume{4} (\bissue{1}): \bfpage{630}--\blpage{655}.
doi:\doiurl{10.3390/polym4010630}.
\burl{https://www.mdpi.com/2073-4360/4/1/630}.
\end{barticle}
\endbibitem

\bibitem[\protect\citeauthoryear{Greenspan}{1976}]{greenspan1976growth}
\begin{barticle}
\bauthor{\bsnm{Greenspan}, \binits{HP}}.
\byear{1976}.
\batitle{On the growth and stability of cell cultures and solid tumors}.
\bjtitle{Journal of theoretical biology}
\bvolume{56} (\bissue{1}): \bfpage{229}--\blpage{242}.
\end{barticle}
\endbibitem

\bibitem[\protect\citeauthoryear{Hao et~al.}{2014}]{hao2014high}
\begin{barticle}
\bauthor{\bsnm{Hao}, \binits{Sijia}}, \bauthor{\binits{Alex~H} \bsnm{Barnett}},
  \bauthor{\binits{Per-Gunnar} \bsnm{Martinsson}}, and \bauthor{\binits{P}
  \bsnm{Young}}.
\byear{2014}.
\batitle{High-order accurate methods for nystr{\"o}m discretization of integral
  equations on smooth curves in the plane}.
\bjtitle{Advances in Computational Mathematics}
\bvolume{40} (\bissue{1}): \bfpage{245}--\blpage{272}.
\end{barticle}
\endbibitem

\bibitem[\protect\citeauthoryear{Hao et~al.}{2018}]{Wenrui2018}
\begin{barticle}
\bauthor{\bsnm{Hao}, \binits{Wenrui}}, \bauthor{\binits{Bei} \bsnm{Hu}},
  \bauthor{\binits{Shuwang} \bsnm{Li}}, and \bauthor{\binits{Lingyu}
  \bsnm{Song}}.
\byear{2018}.
\batitle{Convergence of boundary integral method for a free boundary system}.
\bjtitle{Journal of Computational and Applied Mathematics}
\bvolume{334}: \bfpage{128}--\blpage{157}.
doi:\doiurl{https://doi.org/10.1016/j.cam.2017.11.016}.
\end{barticle}
\endbibitem

\bibitem[\protect\citeauthoryear{He et~al.}{2012}]{he2012modeling}
\begin{barticle}
\bauthor{\bsnm{He}, \binits{Andong}}, \bauthor{\binits{John} \bsnm{Lowengrub}},
  and \bauthor{\binits{Andrew} \bsnm{Belmonte}}.
\byear{2012}.
\batitle{Modeling an elastic fingering instability in a reactive hele-shaw
  flow}.
\bjtitle{SIAM Journal on Applied Mathematics}
\bvolume{72} (\bissue{3}): \bfpage{842}--\blpage{856}.
\end{barticle}
\endbibitem

\bibitem[\protect\citeauthoryear{Helfrich}{1973}]{Helfrich}
\begin{botherref}
\oauthor{\bsnm{Helfrich}, \binits{W.}}
1973.
Elastic properties of lipid bilayers: theory and possible experiments.
\textit{Z Naturf C.}.
\end{botherref}
\endbibitem

\bibitem[\protect\citeauthoryear{Hou et~al.}{1994}]{hou1994removing}
\begin{barticle}
\bauthor{\bsnm{Hou}, \binits{Thomas~Y}}, \bauthor{\binits{John~S}
  \bsnm{Lowengrub}}, and \bauthor{\binits{Michael~J} \bsnm{Shelley}}.
\byear{1994}.
\batitle{Removing the stiffness from interfacial flows with surface tension}.
\bjtitle{Journal of Computational Physics}
\bvolume{114} (\bissue{2}): \bfpage{312}--\blpage{338}.
\end{barticle}
\endbibitem

\bibitem[\protect\citeauthoryear{Hou et~al.}{2001}]{hou2001boundary}
\begin{barticle}
\bauthor{\bsnm{Hou}, \binits{Thomas~Y}}, \bauthor{\binits{John~S}
  \bsnm{Lowengrub}}, and \bauthor{\binits{Michael~J} \bsnm{Shelley}}.
\byear{2001}.
\batitle{Boundary integral methods for multicomponent fluids and multiphase
  materials}.
\bjtitle{Journal of Computational Physics}
\bvolume{169} (\bissue{2}): \bfpage{302}--\blpage{362}.
\end{barticle}
\endbibitem

\bibitem[\protect\citeauthoryear{Krasny}{1986}]{krasny1986study}
\begin{barticle}
\bauthor{\bsnm{Krasny}, \binits{Robert}}.
\byear{1986}.
\batitle{A study of singularity formation in a vortex sheet by the point-vortex
  approximation}.
\bjtitle{Journal of Fluid Mechanics}
\bvolume{167}: \bfpage{65}--\blpage{93}.
\end{barticle}
\endbibitem

\bibitem[\protect\citeauthoryear{Kress}{1995}]{kress1995numerical}
\begin{barticle}
\bauthor{\bsnm{Kress}, \binits{Rainer}}.
\byear{1995}.
\batitle{On the numerical solution of a hypersingular integral equation in
  scattering theory}.
\bjtitle{Journal of computational and applied mathematics}
\bvolume{61} (\bissue{3}): \bfpage{345}--\blpage{360}.
\end{barticle}
\endbibitem

\bibitem[\protect\citeauthoryear{Kress}{2013}]{kress2013linear}
\begin{bbook}
\bauthor{\bsnm{Kress}, \binits{Rainer}}.
\byear{2013}.
\bbtitle{Linear integral equations},
Vol. \bseriesno{82}.
\bpublisher{Springer}.
\end{bbook}
\endbibitem

\bibitem[\protect\citeauthoryear{Leo et~al.}{2000}]{leo2000microstructural}
\begin{barticle}
\bauthor{\bsnm{Leo}, \binits{Perry~H}}, \bauthor{\binits{John~S}
  \bsnm{Lowengrub}}, and \bauthor{\binits{Qing} \bsnm{Nie}}.
\byear{2000}.
\batitle{Microstructural evolution in orthotropic elastic media}.
\bjtitle{Journal of Computational Physics}
\bvolume{157} (\bissue{1}): \bfpage{44}--\blpage{88}.
\end{barticle}
\endbibitem

\bibitem[\protect\citeauthoryear{Li and Li}{2011}]{li2011boundary}
\begin{barticle}
\bauthor{\bsnm{Li}, \binits{Shuwang}}, and \bauthor{\binits{Xiaofan}
  \bsnm{Li}}.
\byear{2011}.
\batitle{A boundary integral method for computing the dynamics of an epitaxial
  island}.
\bjtitle{SIAM Journal on Scientific Computing}
\bvolume{33} (\bissue{6}): \bfpage{3282}--\blpage{3302}.
\end{barticle}
\endbibitem

\bibitem[\protect\citeauthoryear{Li et~al.}{2007}]{ShuwangJCP}
\begin{barticle}
\bauthor{\bsnm{Li}, \binits{Shuwang}}, \bauthor{\binits{John~S.}
  \bsnm{Lowengrub}}, and \bauthor{\binits{Perry~H.} \bsnm{Leo}}.
\byear{2007}.
\batitle{A rescaling scheme with application to the long-time simulation of
  viscous fingering in a hele--shaw cell}.
\bjtitle{Journal of Computational Physics}
\bvolume{225} (\bissue{1}): \bfpage{554}--\blpage{567}.
\end{barticle}
\endbibitem

\bibitem[\protect\citeauthoryear{Liu and Li}{2014}]{Kai2014}
\begin{barticle}
\bauthor{\bsnm{Liu}, \binits{Kai}}, and \bauthor{\binits{Shuwang} \bsnm{Li}}.
\byear{2014}.
\batitle{Nonlinear simulation of vesicle wrinkling}.
\bjtitle{Mathematical Methods in the Applied Sciences}
\bvolume{8}: \bfpage{1093}--\blpage{1112}.
\end{barticle}
\endbibitem

\bibitem[\protect\citeauthoryear{Lodish et~al.}{2008}]{lodish2008molecular}
\begin{bbook}
\bauthor{\bsnm{Lodish}, \binits{Harvey}}, \bauthor{\binits{James~E}
  \bsnm{Darnell}}, \bauthor{\binits{Arnold} \bsnm{Berk}},
  \bauthor{\binits{Chris~A} \bsnm{Kaiser}}, \bauthor{\binits{Monty}
  \bsnm{Krieger}}, \bauthor{\binits{Matthew~P} \bsnm{Scott}},
  \bauthor{\binits{Anthony} \bsnm{Bretscher}}, \bauthor{\binits{Hidde}
  \bsnm{Ploegh}}, \bauthor{\binits{Paul} \bsnm{Matsudaira}}, \betal.
\byear{2008}.
\bbtitle{Molecular cell biology}.
\bpublisher{Macmillan}.
\end{bbook}
\endbibitem

\bibitem[\protect\citeauthoryear{Lowengrub
  et~al.}{2009}]{lowengrub2009nonlinear}
\begin{barticle}
\bauthor{\bsnm{Lowengrub}, \binits{John~S}}, \bauthor{\binits{Hermann~B}
  \bsnm{Frieboes}}, \bauthor{\binits{Fang} \bsnm{Jin}},
  \bauthor{\binits{Yao-Li} \bsnm{Chuang}}, \bauthor{\binits{Xiaolong}
  \bsnm{Li}}, \bauthor{\binits{Paul} \bsnm{Macklin}},
  \bauthor{\binits{Steven~M} \bsnm{Wise}}, and \bauthor{\binits{Vittorio}
  \bsnm{Cristini}}.
\byear{2009}.
\batitle{Nonlinear modelling of cancer: bridging the gap between cells and
  tumours}.
\bjtitle{Nonlinearity}
\bvolume{23} (\bissue{1}): \bfpage{1}.
\end{barticle}
\endbibitem

\bibitem[\protect\citeauthoryear{Macklin and
  Lowengrub}{2007}]{macklin2007nonlinear}
\begin{barticle}
\bauthor{\bsnm{Macklin}, \binits{Paul}}, and \bauthor{\binits{John}
  \bsnm{Lowengrub}}.
\byear{2007}.
\batitle{Nonlinear simulation of the effect of microenvironment on tumor
  growth}.
\bjtitle{Journal of theoretical biology}
\bvolume{245} (\bissue{4}): \bfpage{677}--\blpage{704}.
\end{barticle}
\endbibitem

\bibitem[\protect\citeauthoryear{Martensen}{1963}]{martensen1963methode}
\begin{barticle}
\bauthor{\bsnm{Martensen}, \binits{Erich}}.
\byear{1963}.
\batitle{{\"U}ber eine methode zum r{\"a}umlichen neumannschen problem mit
  einer anwendung f{\"u}r torusartige berandungen}.
\bjtitle{Acta mathematica}
\bvolume{109} (\bissue{1}): \bfpage{75}--\blpage{135}.
\end{barticle}
\endbibitem

\bibitem[\protect\citeauthoryear{Mason~B.N.}{2012}]{Mason2012}
\begin{botherref}
\oauthor{\bsnm{Mason B.~N.} \bsuffix{Califano J.~P.}, \binits{Reinhart-King
  C.~A.}}
2012.
Matrix stiffness: A regulator of cellular behavior and tissue formation.
\textit{Bhatia S. (eds) Engineering Biomaterials for Regenerative Medicine}.
\end{botherref}
\endbibitem

\bibitem[\protect\citeauthoryear{Mikucki and Zhou}{2017}]{Zhou20171}
\begin{barticle}
\bauthor{\bsnm{Mikucki}, \binits{M.}}, and \bauthor{\binits{Y.} \bsnm{Zhou}}.
\byear{2017}.
\batitle{Curvature-driven molecular flow on membrane surface}.
\bjtitle{SIAM Journal on Applied Mathematics}
\bvolume{77} (\bissue{5}): \bfpage{1587}--\blpage{1605}.
doi:\doiurl{10.1137/16M1076551}.
\end{barticle}
\endbibitem

\bibitem[\protect\citeauthoryear{Perthame et~al.}{2014}]{perthame2014hele}
\begin{barticle}
\bauthor{\bsnm{Perthame}, \binits{Beno{\^\i}t}}, \bauthor{\binits{Fernando}
  \bsnm{Quir{\'o}s}}, and \bauthor{\binits{Juan~Luis} \bsnm{V{\'a}zquez}}.
\byear{2014}.
\batitle{The hele--shaw asymptotics for mechanical models of tumor growth}.
\bjtitle{Archive for Rational Mechanics and Analysis}
\bvolume{212} (\bissue{1}): \bfpage{93}--\blpage{127}.
\end{barticle}
\endbibitem

\bibitem[\protect\citeauthoryear{Pham et~al.}{2010}]{pham2010predictions}
\begin{barticle}
\bauthor{\bsnm{Pham}, \binits{Kara}}, \bauthor{\binits{Hermann~B}
  \bsnm{Frieboes}}, \bauthor{\binits{Vittorio} \bsnm{Cristini}}, and
  \bauthor{\binits{John} \bsnm{Lowengrub}}.
\byear{2010}.
\batitle{Predictions of tumour morphological stability and evaluation against
  experimental observations}.
\bjtitle{Journal of the Royal Society Interface}
\bvolume{8} (\bissue{54}): \bfpage{16}--\blpage{29}.
\end{barticle}
\endbibitem

\bibitem[\protect\citeauthoryear{Pham et~al.}{2018}]{pham2018nonlinear}
\begin{botherref}
\oauthor{\bsnm{Pham}, \binits{Kara}}, \oauthor{\binits{Emma} \bsnm{Turian}},
  \oauthor{\binits{Kai} \bsnm{Liu}}, \oauthor{\binits{Shuwang} \bsnm{Li}}, and
  \oauthor{\binits{John} \bsnm{Lowengrub}}.
2018.
Nonlinear studies of tumor morphological stability using a two-fluid flow
  model.
\textit{Journal of mathematical biology}.
\end{botherref}
\endbibitem

\bibitem[\protect\citeauthoryear{Pozrikidis}{1992}]{pozrikidis1992boundary}
\begin{bbook}
\bauthor{\bsnm{Pozrikidis}, \binits{Constantine}}.
\byear{1992}.
\bbtitle{Boundary integral and singularity methods for linearized viscous
  flow}.
\bpublisher{Cambridge University Press}.
\end{bbook}
\endbibitem

\bibitem[\protect\citeauthoryear{Roose et~al.}{2007}]{roose2007mathematical}
\begin{barticle}
\bauthor{\bsnm{Roose}, \binits{Tiina}}, \bauthor{\binits{S~Jonathan}
  \bsnm{Chapman}}, and \bauthor{\binits{Philip~K} \bsnm{Maini}}.
\byear{2007}.
\batitle{Mathematical models of avascular tumor growth}.
\bjtitle{SIAM review}
\bvolume{49} (\bissue{2}): \bfpage{179}--\blpage{208}.
\end{barticle}
\endbibitem

\bibitem[\protect\citeauthoryear{Saad and Schultz}{1986}]{saad1986gmres}
\begin{barticle}
\bauthor{\bsnm{Saad}, \binits{Youcef}}, and \bauthor{\binits{Martin~H}
  \bsnm{Schultz}}.
\byear{1986}.
\batitle{Gmres: A generalized minimal residual algorithm for solving
  nonsymmetric linear systems}.
\bjtitle{SIAM Journal on scientific and statistical computing}
\bvolume{7} (\bissue{3}): \bfpage{856}--\blpage{869}.
\end{barticle}
\endbibitem

\bibitem[\protect\citeauthoryear{Sekerka}{1968}]{Sekerka1968}
\begin{barticle}
\bauthor{\bsnm{Sekerka}, \binits{R.~F.}}
\byear{1968}.
\batitle{Morphological stability}.
\bjtitle{Journal of Crystal Growth}
\bvolume{3}: \bfpage{71}--\blpage{81}.
\end{barticle}
\endbibitem

\bibitem[\protect\citeauthoryear{Sohn et~al.}{2012}]{Sohn2012AxisymmetricMV}
\begin{barticle}
\bauthor{\bsnm{Sohn}, \binits{Jinsun}}, \bauthor{\binits{Shuwang} \bsnm{Li}},
  \bauthor{\binits{Xiaofan} \bsnm{Li}}, and \bauthor{\binits{John~S.}
  \bsnm{Lowengrub}}.
\byear{2012}.
\batitle{Axisymmetric multicomponent vesicles: A comparison of hydrodynamic and
  geometric models.}
\bjtitle{International journal for numerical methods in biomedical engineering}
\bvolume{28}: \bfpage{346}--\blpage{68}.
\end{barticle}
\endbibitem

\bibitem[\protect\citeauthoryear{Turian et~al.}{2018}]{turian2018morphological}
\begin{botherref}
\oauthor{\bsnm{Turian}, \binits{Emma}}, \oauthor{\binits{Kai} \bsnm{Liu}},
  \oauthor{\binits{John} \bsnm{Lowengrub}}, and \oauthor{\binits{Shuwang}
  \bsnm{Li}}.
2018.
Morphological stability of an elastic tumor--host interface.
\textit{Journal of Computational and Applied Mathematics}.
\end{botherref}
\endbibitem

\bibitem[\protect\citeauthoryear{Wei}{2010}]{Wei2010}
\begin{botherref}
\oauthor{\bsnm{Wei}, \binits{GW.}}
2010.
Differential geometry based multiscale models.
\textit{Bull. Math. Biol.}.
\end{botherref}
\endbibitem

\bibitem[\protect\citeauthoryear{Zhao et~al.}{2016}]{zhao2016nonlinear}
\begin{barticle}
\bauthor{\bsnm{Zhao}, \binits{Meng}}, \bauthor{\binits{Andrew}
  \bsnm{Belmonte}}, \bauthor{\binits{Shuwang} \bsnm{Li}},
  \bauthor{\binits{Xiaofan} \bsnm{Li}}, and \bauthor{\binits{John}
  \bsnm{Lowengrub}}.
\byear{2016}.
\batitle{Nonlinear simulations of elastic fingering in a hele-shaw cell}.
\bjtitle{Journal of Computational and Applied Mathematics}
\bvolume{307}: \bfpage{394}--\blpage{407}.
\end{barticle}
\endbibitem

\bibitem[\protect\citeauthoryear{Zhao et~al.}{2017}]{Meng2017}
\begin{barticle}
\bauthor{\bsnm{Zhao}, \binits{Meng}}, \bauthor{\binits{Wenjun} \bsnm{Ying}},
  \bauthor{\binits{John} \bsnm{Lowengrub}}, and \bauthor{\binits{Shuwang}
  \bsnm{Li}}.
\byear{2017}.
\batitle{An efficient adaptive rescaling scheme for computing moving interface
  problems}.
\bjtitle{Communications in Computational Physics}
\bvolume{21} (\bissue{3}): \bfpage{679}--\blpage{691}.
doi:\doiurl{10.4208/cicp.OA-2016-0040}.
\end{barticle}
\endbibitem

\bibitem[\protect\citeauthoryear{Zhu et~al.}{1996}]{zhu1996efficient}
\begin{barticle}
\bauthor{\bsnm{Zhu}, \binits{Jingyi}}, \bauthor{\binits{Xinfu} \bsnm{Chen}},
  and \bauthor{\binits{Thomas~Y} \bsnm{Hou}}.
\byear{1996}.
\batitle{An efficient boundary integral method for the mullins--sekerka
  problem}.
\bjtitle{Journal of Computational Physics}
\bvolume{127} (\bissue{2}): \bfpage{246}--\blpage{267}.
\end{barticle}
\endbibitem

\end{thebibliography}
\nocite{*}

\end{document}